\documentclass{amsart}
\usepackage{amscd,amssymb,hyperref}

\theoremstyle{plain}
\newtheorem{prop}[subsection]{Proposition}
\newtheorem{thm}[subsection]{Theorem}
\newtheorem{lem}[subsection]{Lemma}
\newtheorem{cor}[subsection]{Corollary}
\newtheorem*{theo}{Theorem}

\theoremstyle{remark}
\newtheorem{rem}[subsection]{Remark}

\theoremstyle{definition}
\newtheorem{exm}[subsection]{Example}

\numberwithin{equation}{section}
\renewcommand{\labelenumi}{{\bfseries\arabic{enumi}.}}

\renewcommand{\b}[1]{\mathbf{#1}}
\newcommand{\smallmath}[1]{\text{{\footnotesize{$#1$}}}}

\newcommand{\A}{{\mathcal A}}
\newcommand{\Ai}{{\mathcal A}_\infty}
\newcommand{\B}{{\mathcal B}}
\newcommand{\FF}{{\mathcal F}}
\newcommand{\cE}{{\mathcal E}}
\newcommand{\cG}{{\mathcal G}}
\newcommand{\LL}{{\mathcal L}}
\newcommand{\cS}{{\mathcal S}}
\newcommand{\cT}{{\mathcal T}}

\newcommand{\Z}{{\mathbb Z}}
\newcommand{\C}{{\mathbb C}}
\newcommand{\CP}{{\mathbb{CP}}}
\newcommand{\T}{{({\mathbb C}^*)^n}}
\newcommand{\bS}{{\mathbb S}}

\newcommand{\sfB}{{\sf B}}
\newcommand{\sfD}{{\sf D}}
\newcommand{\sfE}{{\sf E}}

\newcommand{\sfM}{{\sf M}}
\newcommand{\sfY}{{\sf Y}}

\newcommand{\D}{{\Delta}}
\newcommand{\la}{{\lambda }}
\newcommand{\bul}{{\bullet }}

\renewcommand{\c}{{\gamma }}
\renewcommand{\ll}{{\ell }}

\DeclareMathOperator{\ind}{ind}
\DeclareMathOperator{\dep}{dep}
\DeclareMathOperator{\ii}{i}
\DeclareMathOperator{\id}{id}

\DeclareMathOperator{\Aut}{Aut}
\DeclareMathOperator{\End}{End}
\DeclareMathOperator{\Hom}{Hom}

\begin{document}

\title[Gauss-Manin connections for arrangements, II]
{Gauss-Manin connections for arrangements, II \\ Nonresonant weights}
\author[D.~Cohen]{Daniel C.~Cohen$^\dag$}
\address{Department of Mathematics, Louisiana State University,
Baton Rouge, LA 70803}
\email{\href{mailto:cohen@math.lsu.edu}{cohen@math.lsu.edu}}
\urladdr{\href{http://www.math.lsu.edu/~cohen/}
{http://www.math.lsu.edu/\~{}cohen}}
\thanks{{$^\dag$}Partially supported by Louisiana Board of Regents grant
LEQSF(1999-2002)-RD-A-01
and by National Security Agency grant MDA904-00-1-0038}

\author[P.~Orlik]{Peter Orlik$^\ddag$}
\address{Department of Mathematics, University of Wisconsin,
Madison, WI 53706}
\email{\href{mailto:orlik@math.wisc.edu}{orlik@math.wisc.edu}}
\thanks{{$^\ddag$}Partially supported by National Security
Agency grant MDA904-02-1-0019}

\subjclass[2000]{32S22, 14D05, 52C35, 55N25}

\keywords{hyperplane arrangement, local system, Gauss-Manin
connection}

\begin{abstract}
We study the Gauss-Manin connection for the moduli space of an
arrangement of complex hyperplanes in the cohomology of a nonresonant
complex rank one local system.  Aomoto and Kita determined this
Gauss-Manin connection for arrangements in general position.  We use
their results and an algorithm constructed in this paper to determine
this Gauss-Manin connection for all arrangements.
\end{abstract}

\date{July 4, 2002}

\maketitle

\section{Introduction}
\label{sec:intro}
Let $\A=\{H_1,\dots,H_n\}$ be an arrangement of $n$ ordered
hyperplanes in $\C^\ll$, and let $\LL$ be a local system of
coefficients on $\sfM=M(\A)=\C^\ll\setminus\bigcup_{j=1}^n H_j$, the
complement of $\A$.  The need to calculate the local system cohomology
$H^*(\sfM;\LL)$ arises in various contexts.  For instance, local
systems may be used to study the Milnor fiber of the non-isolated
hypersurface singularity at the origin obtained by coning the
arrangement, see \cite{CS1, CO2}.  In mathematical physics, local
systems on complements of arrangements arise in the Aomoto-Gelfand
theory of multivariable hypergeometric integrals \cite{AK,Gel1,OT2}
and the representation theory of Lie algebras and quantum groups.
These considerations lead to solutions of the Knizhnik-Zamolodchikov
differential equation from conformal field theory, see \cite{SV,Va}.
Here a central problem is the determination of the Gauss-Manin
connection on $H^*(\sfM;\LL)$ for certain arrangements, and certain
local systems.

A complex rank one local system on $\sfM$ is determined by a
collection of weights $\la=(\la_1,\dots,\la_n)\in\C^n$.  Associated
to $\la$, we have a representation $\rho:\pi_1(\sfM)\to\C^*$, given by
$\c_j\mapsto \exp(-2\pi\ii\la_j)$ for any meridian loop $\c_j$ about
the hyperplane $H_j$ of $\A$, and an associated rank one local system
$\LL$ on $\sfM$.  Several authors have considered Gauss-Manin
connections with various conditions on the weights: Aomoto~\cite{Ao},
Schechtman and Varchenko \cite{SV,Va}, and Kaneko~\cite{JK} studied
discriminantal arrangements; Kanarek \cite{HK} studied the connection
arising when a single hyperplane in the arrangement is allowed to
move; Aomoto and Kita \cite{AK} determined the Gauss-Manin connection
matrices for general position arrangements; and Terao \cite{T1}
computed these connection matrices for a larger class of arrangements.
In this paper, we determine Gauss-Manin connection matrices for {\em
all} arrangements and all weights $\la$ which satisfy the nonresonance
condition (STV) of \cite{STV} stated in Theorem \ref{thm:vanish}.

Given weights $\la$, we used stratified Morse theory in \cite{C1,CO1}
to construct a complex $K^{\bul}(\A)$ which computes
$H^{*}(\sfM;\LL)$.  This construction is reviewed in
Section~\ref{sec:local}.  In Section \ref{sec:moduli}, we compactify
$\A$ by adding the infinite hyperplane $H_{n+1}$ to obtain the
projective closure $\Ai \subseteq \CP^\ll$.  Two arrangements are said
to have the same combinatorial type if there is an order-preserving
isomorphism between their intersection posets.  If $\cT$ is the
combinatorial type of $\A$, we denote by $\ind(\cT)$ and $\dep(\cT)$
the independence and dependence sets of $\cT$.  These consist of
all $\ll+1$ element subsets of $[n+1]=\{1,2,\dots,n+1\}$ for which the
corresponding hyperplanes in $\Ai$ are independent, and dependent,
respectively.

The sets $\ind(\cT$) and $\dep(\cT)$ may be used to describe the
moduli space of all arrangements of type $\cT$, defined and
investigated by Terao in \cite{T1}.  Let $\sfB(\cT)$ be a smooth,
connected component of this moduli space.  There is a fiber bundle
$\pi:\sfM(\cT) \to \sfB(\cT)$.  The fibers of this bundle,
$\pi^{-1}(\b{x})=\sfM_{\b{x}}$, are complements of arrangements
$\A_{\b{x}}$ of type $\cT$, so they are diffeomorphic to $\sfM$ since
$\sfB(\cT)$ is connected.  Since the bundle $\pi:\sfM(\cT) \to
\sfB(\cT)$ is locally trivial, there is an associated vector bundle
$\pi^{q}:\b{H}^{q}\to\sfB(\cT)$, with fiber
$(\pi^{q})^{-1}(\b{x})=H^{q}(\sfM_{\b{x}};\LL_{\b{x}})$ at
$\b{x}\in\sfB(\cT)$ for each $q$, $0\le q \le \ell$.  The transition
functions of this vector bundle are locally constant.  Fixing a
basepoint $\b{x}\in\sfB(\cT)$, the operation of parallel translation
of fibers over curves in $\sfB(\cT)$ in the vector bundle
$\pi^{q}:\b{H}^{q}\to\sfB(\cT)$ provides a complex representation
\begin{equation} \label{eq:irep}
\Psi^{q}_{\cT}:\pi_{1}(\sfB(\cT),\b{x}) \longrightarrow
\Aut_\C(H^{q}(\sfM_{\b{x}};\LL_{\b{x}})).
\end{equation}

A key idea in this paper is to extract information
about arbitrary arrangements using information about
general position arrangements. Let  $\cG$ denote the
combinatorial type of a general position arrangement
of $n$ hyperplanes in $\C^\ll$, so $\dep(\cG)=\emptyset$.
Let $\sfB(\cG)$ be its moduli space, $\b{y}\in \sfB(\cG)$,
$\B^{n,\ll}_\b{y}$ the corresponding arrangement,
with complement
$\sfE^{n,\ll}_\b{y}=M(\B^{n,\ll}_\b{y})$.
In Section \ref{sec:reps}, we compare the representation
\eqref{eq:irep} with the
analogous representation $\Psi_\cG^q:\pi_1(\sfB(\cG),\b{y}) \to
\Aut_\C(H^{q}(\sfE^{n,\ll}_{\b{y}};\LL_{\b{y}}))$.  To this
end, we introduce a space $\sfY(\cT)$ which contains $\sfB(\cT)$ and
$\sfB(\cG)$ as disjoint subspaces.  Let $i_{\cT}:\sfB(\cG) \to
\sfY(\cT)$ and $j_{\cT}:\sfB(\cT) \to \sfY(\cT)$ denote the natural
inclusions.  Given $\c\in\pi_1(\sfB(\cT),\b{x})$, we show that there
is a {\em compatible} $\Gamma \in \pi_1(\sfB(\cG),\b{y})$, a loop for
which the homology classes $[\c]$ and $[\Gamma]$ satisfy
$(i_\cT)_*([\Gamma])=(j_\cT)_*([\c])$.  For compatible loops, the
automorphisms $\Psi^q_\cT(\c)$ and $\Psi^q_\cG(\Gamma)$
are related by
\begin{equation}\label{eq:icong}
\overline{s}^q_\b{x} \circ \Psi^q_\cG(\Gamma) \sim \Psi^q_\cT(\c) \circ
\overline{s}^q_\b{x}.
\end{equation}
Here $\overline{s}^q_\b{x}:H^q(M(\B^{n,\ll}_\b{y});\LL_\b{y}) \to
H^q(M(\A_\b{x});\LL_\b{x})$ is the map in local system cohomology
induced by a certain chain map $s^\bul_\b{x}:K^\bul(\B^{n,\ll}_\b{y})
\to K^\bul(\A_\b{x})$, and $\sim$ denotes equality up to conjugacy.

These compatible classes are studied further in Section
\ref{sec:compatible}. The moduli space $\sfB(\cG)$ of general position
arrangements is the complement of a divisor $\sfD(\cG)= \bigcup_J
\sfD_J$ in $(\CP^\ll)^n$.  The components $\sfD_J$ of this divisor are
irreducible hypersurfaces indexed by $\ll +1$ element subsets $J$ of
$[n+1]$.  For each such $J$, let $\Gamma_{\!J}$ be a meridian loop in
$\sfB(\cG)$ about a generic point in $\sfD_J$.  We prove that the
homology group $H_1(\sfB(\cG))$ is generated by the classes
$[\Gamma_{\!J}]$.   Note that $\Gamma_{\!J}$ is a loop in $\sfY(\cT)$ for
any
combinatorial type $\cT$ since $\sfB(\cG) \subseteq \sfY(\cT)$.  We
say that the type $\cT$ covers $\cT'$ if $\sfB(\cT')$ has complex
codimension one in the closure $\overline{\sfB}(\cT)$.
If $\c\in\pi_1(\sfB(\cT),\b{x})$ is a
simple loop in $\sfB(\cT)$ about a generic point in $\sfB(\cT')$, we
show that the homology class $[\c]$ satisfies
\begin{equation} \label{eq:icompat}
(j_\cT)_*([\c]) = \sum_{J\in \dep(\cT',\cT)} m_J
\cdot [\Gamma_{\!J}],
\end{equation}
where $\dep(\cT',\cT)=\dep(\cT') \setminus \dep(\cT)$, and $m_J$ is
the order of vanishing of the restriction of a defining polynomial for
$\sfD_J$ to $\overline{\sfB}(\cT)$ along $\sfB(\cT')$.

If $\LL$ is a nontrivial rank one local system on the complement
$\sfE^{n,\ll}$ of a general position arrangement of $n$ hyperplanes in
$\C^\ll$, then the cohomology $H^q(\sfE^{n,\ll};\LL)$ vanishes in all
dimensions $q$ except possibly $q=\ll$.  For an arbitrary arrangement
$\A$, a collection of weights $\la$ and the corresponding local system
$\LL$ are called {\em nonresonant} if the cohomology $H^q(\sfM;\LL)$
vanishes in all dimensions $q$ except possibly $q=\ll$.  We review the
(STV)-condition for $\cT$--nonresonance \cite{STV} and the
$\beta$\textbf{nbc} basis for the nonvanishing group \cite{FT} in
Section \ref{sec:weights}.  We show that for $\cT$--nonresonant
weights $\la$, the chain map $s^\bul:K^\bul(\B^{n,\ll}) \to
K^\bul(\A)$  induces a surjection
$\overline{s}^\ll:H^\ll(\sfE^{n,\ll},\LL) \twoheadrightarrow
H^\ll(\sfM;\LL)$.

The vector bundle $\pi^q:\b{H}^q\to\sfB(\cT)$, with fiber
$(\pi^q)^{-1}(\b{x})=H^q(\sfM_\b{x};\LL_\b{x})$ at
$\b{x}\in\sfB(\cT)$, supports a Gauss-Manin connection corresponding
to the representation  \eqref{eq:irep}.  We study
this connection in Section \ref{sec:GM} for $\cT$--nonresonant
weights, where only $q=\ll$ is relevant.
Let $\c\in\pi_1(\sfB(\cT),\b{x})$ be represented by
$g:\bS^1\to\sfB(\cT)$.
A connection matrix, $\Omega_\cT(\c)$, for the induced
connection on the pullback of the vector bundle
$\pi^\ll:\b{H}^\ll\to\sfB(\cT)$ along $g$ is called a
Gauss-Manin connection matrix. Let
$\Gamma\in\pi_1(\sfB(\cG),\b{y})$ be a loop whose homology class is
given by \eqref{eq:icompat} above, $[\Gamma]=\sum_{J\in \dep(\cT',\cT)} m_J
\cdot [\Gamma_{\!J}]$.  Then $\c$ and $\Gamma$ are compatible, and it
follows from \eqref{eq:icong} that
$\overline{s}^\ll_\b{x} \circ \Psi^\ll_\cG(\Gamma) \sim
\Psi^\ll_\cT(\c) \circ \overline{s}^\ll_\b{x}$.
Since the homomorphism
$\overline{s}^\ll_\b{x}:H^\ll(\sfE^{n,\ll}_\b{y},\LL_\b{y}) \to
H^\ll(\sfM_\b{x};\LL_\b{x})$ is surjective for $\cT$--nonresonant
weights, this determines the Gauss-Manin connection matrix
$\Omega_\cT(\c)$.

Aomoto and Kita \cite{AK} obtained explicit formulas for the
Gauss-Manin connection matrices $\Omega_\cG(\Gamma_J)$ in the case of
general position arrangements, see also \cite{OT2}.  We combine our
results with these formulas to obtain the main result of this paper.
It is a combinatorial algorithm for calculating Gauss-Manin connection
matrices for an arbitrary combinatorial type $\cT$, and
$\cT$--nonresonant weights, in terms of those for the type $\cG$ of
general position arrangements.  Let $\cT$ be a combinatorial type
which covers the type $\cT'$.  Let $\c \in \pi_1(\sfB(\cT),\b{x})$ a
simple loop in $\sfB(\cT)$ about a generic point in $\sfB(\cT')$.  We
agree to write $ \Omega_\cT(\cT')$ in place of $ \Omega_\cT(\c)$ in
this situation.  Similarly, we write $ \Omega_\cG(J)$ in place of $
\Omega_\cG(\Gamma_{\!J})$.  Fix $\cT$--nonresonant (and hence
$\cG$--nonresonant) weights $\la$.  For $\b{x} \in \sfB(\cT)$ and
$\b{y} \in \sfB(\cG)$, the corresponding local system cohomology
groups of the fibers, $H^q(\sfM_\b{x};\LL_\b{x})$ and
$H^q(\sfE^{n,\ell}_\b{y};\LL_\b{y})$, depend only on the combinatorial
types $\cT$ and $\cG$ and vanish for $q \neq \ell$.  So we write
$H^\ell(\sfM_\b{x};\LL_\b{x})=H^\ell(\cT)$ and
$H^\ell(\sfE^{n,\ell}_\b{y};\LL_\b{y})=H^\ell(\cG)$.  This notation
makes it clear that the Gauss-Manin connection matrix $\Omega_\cT(\cT')$
depends only on
the combinatorial types involved.
The matrix $\Omega_\cT(\cT')$ is obtained by the
combinatorial formula recorded below.

\begin{theo}
Let $\cT$ be a combinatorial type which covers the type $\cT'$.  Let
$\la$ be a collection of $\cT$--nonresonant weights, and let $P(\cT)$
be the matrix of the surjection $H^\ll(\cG) \twoheadrightarrow
H^\ll(\cT)$ in the respective $\beta$\textbf{nbc} bases.  Then a
Gauss-Manin connection matrix $\Omega_\cT(\cT')$ is determined by the
matrix equation
\[
P(\cT) \cdot \Omega_\cT(\cT') = \Bigl( \sum_{J \in \dep(\cT',\cT)}
m_J \cdot \Omega_\cG(J) \Bigr) \cdot P(\cT).
\]
\end{theo}
\noindent Illustrations of the algorithm provided by this result are given
in
Section \ref{sec:examples}.

\section{Preliminary Results}
\label{sec:local}
In this section, we record a number of results concerning a fixed
arrangement.  We use notation and results from \cite{OT1,OT2}.

Let $\A=\{H_1,\dots,H_n\}$ be an arrangement of $n$ ordered
hyperplanes in $\C^\ll$.  Assume that the first $\ll$ hyperplanes of
$\A$ are linearly independent, so $\A$ is essential.  Choose
coordinates $u_1,\dots,u_\ll$ for $\C^\ll$, and for each hyperplane
$H_j$ of $\A$, let $\alpha_j=\alpha_j(u_1,\dots,u_\ll)$ be a linear
polynomial such that $H_j$ is defined by the vanishing of $\alpha_j$.
Then $Q(\A)=\prod_{j=1}^n \alpha_j(u_1,\dots,u_\ll)$ is a defining
polynomial for $\A$, and the complement $\sfM=M(\A)$ is given by
$\sfM=\C^\ll \setminus \bigcup_{j=1}^n H_j = \C^\ll \setminus
Q(\A)^{-1}(0)$.

Choose coordinates $u_1,\dots,u_n$ for $\C^n$, and consider the
arrangement $\B^n$ in $\C^n$ defined by $Q(\B^n)=\prod_{j=1}^\ll
\alpha_j(u_1,\dots,u_\ll) \prod_{j=\ll+1}^n u_j$.  Since the first
$\ll$ hyperplanes of $\A$ are linearly independent, $\B^n$ is an
arrangement of $n$ linearly independent hyperplanes in $\C^n$, a
Boolean arrangement, and the complement $\sfE^n=M(\B^n)=\T$ is the
complex $n$-torus.  Define $s:\C^\ll \to \C^n$ by
\begin{equation} \label{eq:section}
s(u_1,\dots,u_\ll)=\left(u_1,\dots,u_\ll,\alpha_{\ll+1}(u_1,\dots,u_\ll),
\dots, \alpha_{n}(u_1,\dots,u_\ll)\right).
\end{equation}
Then $s:\sfM \to \sfE^n$, and we have the following.

\begin{prop} \label{prop:slice}
Let $\A$ be an essential arrangement of $n$ hyperplanes in $\C^\ll$.
Then the complement $\sfM$ of $\A$ is an affine section of the complex
$n$-torus $\sfE^n$.
\end{prop}

Let $A(\A)$ be the Orlik-Solomon algebra of $\A$, generated by one
dimensional classes $a_j$, $1\leq j\leq n$.  It is the quotient of the
exterior algebra generated by these classes by a homogeneous ideal,
hence a finite dimensional graded $\C$-algebra.  There is an
isomorphism of graded algebras $H^*(\sfM;\C) \simeq A(\A)$.  In
particular, $\dim A^q(\A)=b_q(\A)$ where $b_q(\A)=\dim H^q(\sfM;\C)$
denotes the $q$-th Betti number of $\sfM$ with trivial local
coefficients $\C$.

This construction is realized topologically by the map $s:\sfM \to
\sfE^n$ defined above.  The cohomology $H^*(\sfE^n;\C)$ is an exterior
algebra, generated by classes $a_j$ dual to meridian loops $\c_j$
about the hyperplanes of $\B^n$ (respectively $\A$), and one can show
that the map $s$ induces a surjection
\[
A(\B^n)\simeq H^*(\sfE^n;\C)\xrightarrow{\ s^*}H^*(\sfM;\C)\simeq
A(\A)
\]
in (de Rham) cohomology.

Let $\B^{n,\ll}$ be a general position arrangement of $n$ hyperplanes
in $\C^\ll$, a generic $\ll$-dimensional section of the Boolean
arrangement $\B^n$ in $\C^n$.  The Orlik-Solomon algebra
$A(\B^{n,\ll})$ is the rank $\ll$ truncation of the exterior algebra
$A(\B^n)$.  We have $A^q(\B^{n,\ll}) = A^q(\B^n)$ for $q \le \ll$ and
$A^q(\B^{n,\ll})=0$ for $q>\ll$.  The complement of the arrangement
$\A$ in $\C^\ll$ has the homotopy type of an $\ll$-dimensional
complex, so the cohomology of the complement vanishes in higher
dimensions, $H^q(\sfM;\C)=A^q(\A)=0$ for $q>\ll$.  These facts yield
the following result.

\begin{prop} \label{prop:GPsurjection}
The map $s:\sfM \to \sfE^n$ induces a surjection
$s^*:A(\B^{n,\ll}) \twoheadrightarrow A(\A)$.
\end{prop}

A collection of weights $\la=(\la_1,\dots,\la_n)\in\C^n$ gives rise to
compatible rank one local systems on the complements $\sfE^n$,
$\sfE^{n,\ll}$, and $\sfM$ of the arrangements $\B^n$, $\B^{n,\ll}$,
and $\A$ as follows.  Associated to $\la$, there is a rank one
representation $\rho:\pi_1(\sfE^n) \to \C^*$ given by $\c_j \mapsto
\exp(-2\pi\ii\la_j)$ for any meridian loop $\c_j$ about the $j$-th
hyperplane of $\B^n$, and a corresponding rank one local system $\LL$
on $\sfE^n$.  Let $s:\sfM \to \sfE^n$ and $\sigma:\sfE^{n,\ll} \to
\sfE^n$ be the maps from Proposition \ref{prop:slice} realizing $\sfM$
and $\sfE^{n,\ll}$ as sections of $\sfE^n$.  Via these maps, there are
induced local systems on $\sfM$ and $\sfE^{n,\ll}$ corresponding to
representations of $\pi_1(\sfM)$ and $\pi_1(\sfE^{n,\ll})$ defined as
above.  For brevity, we also use $\LL$ to denote these local systems.

One can use stratified Morse theory to construct a cochain complex
$K^{\bul}(\A)$, the cohomology of which is naturally isomorphic to the
local system cohomology $H^{*}(\sfM;\LL)$, see \cite{C1,CO1}.
Briefly, let $\emptyset = \FF^{-1} \subset \FF^0 \subset \FF^1 \subset
\dots \subset \FF^\ll = \C^\ll$ be a flag in $\C^\ll$ which is
transverse to the arrangement $\A$, and let $\sfM^q = \FF^q \cap \sfM$
for each~$q$.  Let $K^q=K^q(\A)=H^q(\sfM^q,\sfM^{q-1};\LL)$, and let
$\D^q:H^{q}(\sfM^{q},\sfM^{q-1};\LL) \to
H^{q+1}(\sfM^{q+1},\sfM^q;\LL)$ be the boundary homomorphism of the
triple $(\sfM^{q+1},\sfM^q,\sfM^{q-1})$.

\begin{thm}[\cite{C1}]\label{thm:Kdot}\

\begin{enumerate}
\item \label{item:Kdot1}
For each $q$, $0\le q \le \ll$, we have $H^i(\sfM^q,\sfM^{q-1};\LL) =
0$ if $i \neq q$, and $\dim_\C H^q(\sfM^q,\sfM^{q-1};\LL) = b_q(\A)$
is equal to the $q$-th Betti number of $\sfM$ with trivial local
coefficients $\C$.

\item \label{item:Kdot2}
The system of complex vector spaces and linear maps
$(K^\bul,\D^\bul)$,
\[
K^0 \xrightarrow{\ \D^{0}\ } K^1 \xrightarrow{\ \D^1\ } K^2
\xrightarrow{\phantom{\D^{1}}} \cdots \xrightarrow{\phantom{\D^{1}}}
K^{\ll-1} \xrightarrow{\ \D^{\ll-1}\,} K^\ll,
\]
is a complex $(\D^{q+1}\circ\D^q=0)$.  The cohomology of this complex
is naturally isomorphic to $H^*(\sfM;\LL)$, the cohomology of $\sfM$
with coefficients in $\LL$.
\end{enumerate}
\end{thm}

\section{Moduli Spaces}
\label{sec:moduli}
We now consider families of arrangements with a given combinatorial
type.  Fix a pair $(\ell,n)$ with $n\ge \ell\geq 1$.  We consider
(multi-)arrangements of $n$ ordered hyperplanes in $\C^\ell$.  We
assume that the arrangement $\A$ contains $\ell$ linearly independent
hyperplanes and that these are the first $\ell$ in the linear order.
Let $u_1,\ldots, u_\ell$ be the coordinates of $\C^\ell$.  Choose
linear polynomials $\alpha_{i} = x_{i,0} + \sum_{j=1}^{\ell} x_{i,j}
u_{j} \,\,(i = 1,\dots, n)$ which define the hyperplanes of $\A$.
Note that the matrix $(x_{i,j})$, $1\le i,j \le \ell$ is invertible by
assumption.

We embed the arrangement in projective space and add the hyperplane at
infinity as last in the ordering, $H_{n+1}$.  The moduli space of all
arrangements combinatorially equivalent to $\A$ may be viewed as the
set of matrices
\begin{equation} \label{eq:point}
\b{x}=
\begin{pmatrix}
x_{1,0} & x_{1,1} & \cdots & x_{1,\ell}\\
x_{2,0} & x_{2,1} & \cdots & x_{2,\ell}\\
\vdots  & \vdots  & \ddots & \vdots \\
x_{n,0} & x_{n,1} & \cdots & x_{n,\ell}\\
1       & 0       & \cdots & 0
\end{pmatrix}
\end{equation}
whose rows are elements of $\CP^\ell$, and whose
$(\ell+1)\times(\ell+1)$ minors satisfy certain dependency conditions,
see \cite[Prop.~ 9.2.2]{OT2}.  Given an $\ell+1$ element subset $I$ of
$[n+1]:=\{1,\dots,n,n+1\}$, let $\Delta_{I}=\Delta_{I}(\b{x})$ denote
the determinant of the submatrix of $\b{x}$ with rows specified by
$I$.  Corresponding to each such $\b{x}$, we have an arrangement
$\A_\b{x}$, with hyperplanes defined by the first $n$ rows of the
matrix equation $\b{x} \cdot \tilde{\sf u}=0$, where $\tilde{\sf u} =
\begin{pmatrix} 1 & u_1 & \cdots & u_\ell\end{pmatrix}^\top$, that is
combinatorially equivalent to $\A$.  Let $\sfM_\b{x}=M(\A_\b{x})$ be
the complement of $\A_\b{x}$.

Let $\cT$ denote the combinatorial type of the arrangement $\A$, and
let $\sfB(\cT)$ denote a smooth, connected component of the moduli
space of all arrangements of type $\cT$.  Let
\[
\sfM(\cT) = \{ (\b{x},{\sf u}) \in (\CP^\ell)^n \times \C^\ell \mid
\b{x} \in \sfB(\cT) \ \hbox{and}\ {\sf u} \in \sfM_\b{x}\},
\]
and define $\pi_{\cT}:\sfM(\cT) \to \sfB(\cT)$ by
$\pi_{\cT}(\b{x},{\sf u})=\b{x}$.  A result of Randell \cite{Ra}
implies that $\pi_{\cT}:\sfM(\cT) \to \sfB(\cT)$ is a bundle, with
fiber $\pi_{\cT}^{-1}(\b{x}) = \sfM_\b{x}$, the complement of the
arrangement defined by $\b{x} \in \sfB(\cT)$.

We now construct a number of spaces and bundles related to the moduli
space $\sfB(\cT)$.  For any combinatorial type $\cT$, let $\ind(\cT)$
denote the set of all $\ell+1$ element subsets
$I=\{i_1,\dots,i_{\ell+1}\}$ of $[n+1]$ for which $\Delta_{I}\neq 0$
in type $\cT$.  If $\cT$ is realizable, $\ind(\cT)$ is the set of all
subsets $I$ for which $\{H_{i_1},\dots,H_{i_{\ell+1}}\}$ is linearly
independent in the projective closure of an arrangement $\A$ of type
$\cT$.  Similarly, let $\dep(\cT)$ be the set of all $\ell+1$ element
subsets $J$ of $[n+1]$ for which $\Delta_{J}=0$ in type $\cT$.

Given a type $\cT$, let
\[
\sfY(\cT) = \{\b{x} \in (\CP^{\ell})^{n} \mid \Delta_{I}(\b{x}) \neq 0
\text{ for } I \in \ind(\cT)\}.
\]
Then the moduli space $\sfB(\cT)$ may be realized as
\[
\sfB(\cT) = \{\b{x} \in \sfY(\cT) \mid \Delta_{J}(\b{x})=0
\text{ for } J \in \dep(\cT)\}.
\]
Note that if $\cG$ is the type of a general position arrangement, then
$\dep(\cG)=\emptyset$, so $\sfB(\cG)=\sfY(\cG)$.  For any type $\cT$,
the moduli space $\sfB(\cG)$ of all general position arrangements may
be realized as
\[
\sfB(\cG)=\{\b{x} \in \sfY(\cT) \mid \Delta_{J}(\b{x}) \neq 0
\text{ for } J \in \dep(\cT)\}.
\]
If $\cT \neq \cG$, then $\sfB(\cT)$ and $\sfB(\cG)$ are disjoint
subspaces of $\sfY(\cT)$.  Let $i_{\cT}:\sfB(\cG) \to \sfY(\cT)$ and
$j_{\cT}:\sfB(\cT) \to \sfY(\cT)$ denote the natural inclusions.

Let $I_{0}=\{1,\dots,\ell,n+1\}$, and define $\cT_{0}$ by
$\ind(\cT_{0})=\{I_{0}\}$.  Then $\sfB(\cT_{0})$ is the moduli space
of all multi-arrangements with $\{H_{1},\dots,H_{\ell},H_{n+1}\}$
linearly independent.  Given $\b{x} \in \sfB(\cT_{0})$, let
$\B^{n}_{\b{x}}$ denote the Boolean arrangement of $n$ hyperplanes in
$\C^{n}$ with defining polynomial $Q(\B^{n}_{\b{x}})$ given by
\[
Q(\B^{n}_{\b{x}}) =
\prod_{i=1}^{\ell}(x_{i,0}+x_{i,1}u_{1}+\dots+x_{i,\ell}u_{\ell})
\prod_{i=\ell+1}^{n} u_{i},
\]
where $u_{1},\dots,u_{\ell},u_{\ell+1},\dots,u_{n}$ are the
coordinates for $\C^{n}$.  Let $\sfE^n_{\b{x}} = \C^{n} \setminus
Q(\B^{n}_{\b{x}})^{-1}(0)$ be the complement of the arrangement
$\B^{n}_{\b{x}}$.  Note that $\sfE^n_{\b{x}} \cong (\C^{*})^{n}$ is a
complex $n$-torus. Define
\[
\sfE(\cT_{0}) = \{(\b{x},{\sf u}) \in (\CP^{\ell})^{n} \times \C^{n}
\mid \b{x} \in \sfY(\cT_{0}) \text{ and } {\sf u} \in
\sfE^n_{\b{x}}\}.
\]

\begin{prop}
The natural map $p_{\cT_{0}}:\sfE(\cT_{0}) \to \sfY(\cT_{0})$ defined
by $p_{\cT_{0}}^{}(\b{x},{\sf u})=\b{x}$ is a fiber bundle projection,
with fiber the complex $n$-torus
$p_{\cT_{0}}^{-1}(\b{x})=\sfE^n_{\b{x}}$.
\end{prop}
\begin{proof}
This may be established using the Thom Isotopy Lemma by modifying the
argument given by Randell \cite{Ra}.
\end{proof}

For each type $\cT$ for which $I_0 \in \ind(\cT)$, let
$k_\cT:\sfY(\cT) \to \sfY(\cT_0)$ denote the natural inclusion.
Recall that $j_\cT:\sfB(\cT) \to \sfY(\cT)$ denotes the inclusion of
the moduli space of $\cT$ in $\sfY(\cT)$.  By pulling back the above
bundle along the inclusion maps $k_\cT:\sfY(\cT) \to \sfY(\cT_0)$ and
$k_\cT \circ j_\cT:\sfB(\cT) \to \sfY(\cT_0)$, we obtain bundles over
$\sfY(\cT)$ and $\sfB(\cT)$, with fiber the complex $n$-torus.

Denote these bundles by
\[
p_{\cT}:\sfE(\cT) \to \sfY(\cT) \quad \text{and} \quad
p'_{\cT}:\sfE'(\cT) \to \sfB(\cT)
\]
respectively.  The total spaces of these bundles may be realized as
\begin{align*}
\sfE(\cT)&= \{(\b{x},{\sf u}) \in (\CP^{\ell})^{n} \times \C^{n} \mid
\b{x} \in \sfY(\cT) \text{ and } {\sf u} \in \sfE^n_{\b{x}}\}\
\text{and}\\
\sfE'(\cT)&=\{(\b{x},{\sf u}) \in (\CP^{\ell})^{n} \times \C^{n} \mid
\b{x} \in \sfB(\cT) \text{ and } {\sf u} \in \sfE^n_{\b{x}}\}.
\end{align*}

The bundle $p'_{\cT}:\sfE'(\cT) \to \sfB(\cT)$ and the moduli space
bundle $\pi_{\cT}:\sfM(\cT) \to \sfB(\cT)$ are related as follows.

\begin{prop} \label{prop:bundlemap}
There is a bundle map $S:\sfM(\cT) \to \sf E'(\cT)$ covering the
identity map of $\sfB(\cT)$.
\end{prop}
\begin{proof}
Let $\b{x} \in \sfB(\cT)$.  Corresponding to $\b{x}$, we have an
arrangement $\A_{\b{x}}$ of type $\cT$ in $\C^{\ell}$, and a Boolean
arrangement $\B^{n}_{\b{x}}$ in $\C^{n}$.  These arrangements have
defining polynomials
\begin{align*}
Q(\A_{\b{x}})&=
\prod_{i=1}^{\ell}(x_{i,0}+x_{i,1}u_{1}+\dots+x_{i,\ell}u_{\ell})
\prod_{i=\ell+1}^{n}(x_{i,0}+x_{i,1}u_{1}+\dots+x_{i,\ell}u_{\ell})
\ \text{and}\\
Q(\B^{n}_{\b{x}})&=
\prod_{i=1}^{\ell}(x_{i,0}+x_{i,1}u_{1}+\dots+x_{i,\ell}u_{\ell})
\prod_{i=\ell+1}^{n} u_{i},
\end{align*}
where $u_{1},\dots,u_{\ell}$ are coordinates for $\C^{\ell}$ and
$u_{1},\dots,u_{n}$ are coordinates for $\C^{n}$.  It follows that the
map $\C^{\ell} \to \C^{n}$, $(u_{1},\dots,u_{\ell}) \mapsto
(u_{1},\dots,u_{\ell},\alpha_{\ll+1},\dots,\alpha_n)$, defined in
\eqref{eq:section} restricts to a map $s^{}_{\b{x}}:\sfM_{\b{x}} \to
\sfE^n_{\b{x}}$ of the complement of $\A_{\b{x}}$ to the complement of
$\B^{n}_{\b{x}}$, see Proposition \ref{prop:slice}.  Defining
$S:\sfM(\cT) \to \sf E'(\cT)$ by
$S(\b{x},{\sf u})=(\b{x},s^{}_{\b{x}}({\sf u}))$ yields the desired
bundle map.
\end{proof}

\section{Representations}
\label{sec:reps}
Let $\A$ be an essential arrangement of $n$ hyperplanes in $\C^{\ell}$
of combinatorial type $\cT$, and let $\sfB(\cT)$ be a smooth,
connected component of the corresponding moduli space.  For each
$\b{x} \in \sfB(\cT)$, a collection of complex weights
$\la=(\la_{1},\dots,\la_{n})$ determines a rank one local system
$\LL_{\b{x}}$ on the complement $\sfM_{\b{x}}$ of the arrangement
$\A_{\b{x}}$, with monodromy $\exp(-2\pi\ii \la_{j})$ about the
hyperplane $H_{j}\in \A_{\b{x}}$.

Since the fiber bundle $\pi:\sfM(\cT) \to \sfB(\cT)$ is locally
trivial, there is an associated vector bundle
$\pi^{q}:\b{H}^{q}\to\sfB(\cT)$, with fiber
$(\pi^{q})^{-1}(\b{x})=H^{q}(\sfM_{\b{x}};\LL_{\b{x}})$ at
$\b{x}\in\sfB(\cT)$ for each $q$, $0\le q \le \ell$.  The transition
functions of this vector bundle are locally constant.  Fixing a
basepoint $\b{x}\in\sfB(\cT)$, the operation of parallel translation
of fibers over curves in $\sfB(\cT)$ in the vector bundle
$\pi^{q}:\b{H}^{q}\to\sfB(\cT)$ provides a complex representation
\begin{equation} \label{eq:Hqrep}
\Psi^{q}:\pi_{1}(\sfB(\cT),\b{x}) \longrightarrow
\Aut_\C(H^{q}(\sfM_{\b{x}};\LL_{\b{x}})).
\end{equation}
Write $\Psi^q = \Psi^q_\cT$ to indicate the dependence of this
representation on the type $\cT$.

By Theorem \ref{thm:Kdot}, the local system cohomology of $\sfM_\b{x}$
may be computed using the Morse theoretic complex $K^\bul(\A_\b{x})$.
The fundmental group of $\sfB(\cT)$ acts by chain automorphisms on
this complex, see \cite[Cor.~3.2]{CO3}, yielding a representation
\begin{equation} \label{eq:KArep}
\Phi^{\bul}:\pi_{1}(\sfB(\cT),\b{x}) \longrightarrow
\Aut_\C(K^{\bul}(\A_{\b{x}})).
\end{equation}
As above, write $\Phi^\bul = \Phi^\bul_\cT$ to indicate the dependence
of this representation on the combinatorial type $\cT$.

\begin{thm}[\cite{CO3}] \label{thm:inducedrep}
The representation $\Psi^{q}_\cT:\pi_{1}(\sfB(\cT),\b{x}) \to
\Aut_\C(H^{q}(\sfM_{\b{x}};\LL_{\b{x}}))$ is induced by the
representation $\Phi^{\bul}_\cT:\pi_{1}(\sfB(\cT),\b{x}) \to
\Aut_\C(K^{\bul}(\A_{\b{x}}))$.
\end{thm}

The constructions of the previous section provide additional
representations of the fundamental group of $\sfB(\cT)$.  Recall that
$I_{0}=\{1,\dots,\ell,n+1\}$, and assume $I_{0} \in \ind(\cT)$.  Then
there is a bundle $p'_{\cT}:\sfE'(\cT) \to \sfB(\cT)$, with fiber
$\sfE^n_{\b{x}}$, the complement of the Boolean arrangement
$\B^{n}_{\b{x}}$ in $\C^{n}$.  Recall also that the the weights $\la$
give rise to a local system on $\sfE^n_\b{x}$, which we also denote by
$\LL_\b{x}$ since it is compatible with the local system on
$\sfM_\b{x}$.  As above, the fundamental group of $\sfB(\cT)$ acts by
chain automorphisms on the cochain complex $K^{\bul}(\B^{n}_{\b{x}})$
associated with this arrangement, yielding a representation
\begin{equation} \label{eq:KBrep}
\widetilde\Phi^{\bul}_\cT:\pi_{1}(\sfB(\cT),\b{x}) \longrightarrow
\Aut_\C(K^{\bul}(\B^{n}_{\b{x}})).
\end{equation}

By Proposition \ref{prop:bundlemap}, there is a bundle map
$s:\sfM(\cT) \to \sf E'(\cT)$, which restricts to
$s_{\b{x}}:\sfM_{\b{x}} \to \sfE^n_{\b{x}}$ on fibers.  The map
$s_{\b{x}}$ induces a chain map
$s_{\b{x}}^{\bul}:K^{\bul}(\B^{n}_{\b{x}}) \to K^{\bul}(\A_{\b{x}})$,
see \cite[Prop.~2.11]{CO1}, which is compatible with the
representations $\Phi^\bul_\cT$ and $\widetilde\Phi^\bul_\cT$ of
\eqref{eq:KArep} and \eqref{eq:KBrep} in the following sense.

\begin{prop}\label{prop:Acd}
For each $\c\in \pi_{1}(\sfB(\cT),\b{x})$, there is a commutative
diagram
\begin{equation*} \label{eq:Acd}
\begin{CD}
K^{\bul}(\B^{n}_{\b{x}}) @>s_{\b{x}}^{\bul}>> K^{\bul}(\A_{\b{x}})\\
@VV\widetilde\Phi^{\bul}_\cT(\c)V      @VV\Phi^{\bul}_\cT(\c)V \\
K^{\bul}(\B^{n}_{\b{x}}) @>s_{\b{x}}^{\bul}>> K^{\bul}(\A_{\b{x}})
\end{CD}
\end{equation*}
\end{prop}

\begin{rem}\label{rem:gp}
Recall that $\cG$ denotes the combinatorial type of a general position
arrangement of $n$ hyperplanes in $\C^\ell$.  If $\b{y}\in \sfB(\cG)$
is a point in the moduli space of all such arrangements, let
$\B^{n,\ell}_\b{y}$ denote the corresponding arrangement.  The above
discussion, when applied to the arrangement $\B^{n,\ell}_\b{y}$,
yields the following commutative diagram
\begin{equation} \label{eq:Bcd}
\begin{CD}
K^{\bul}(\B^{n}_{\b{y}}) @>\sigma_{\b{y}}^{\bul}>>
K^{\bul}(\B^{n,\ell}_{\b{y}})\\
@VV\widetilde\Phi^{\bul}_\cG(\Gamma)V
@VV\Phi^{\bul}_\cG(\Gamma)V \\
K^{\bul}(\B^{n}_{\b{y}}) @>\sigma_{\b{y}}^{\bul}>>
K^{\bul}(\B^{n,\ell}_{\b{y}})
\end{CD}
\end{equation}
where $\Gamma \in \pi_1(\sfB(\cG),\b{y})$ and $\sigma_\b{y}^\bul:
K^{\bul}(\B^{n}_{\b{y}}) \to K^{\bul}(\B^{n,\ell}_{\b{y}})$ is the
chain map induced by the bundle map $\sigma:\sfM(\cG) \to \sfE'(\cG)$
of Proposition \ref{prop:bundlemap} in this particular case.

The complement $\sfE^{n,\ll}_\b{y}=M(\B^{n,\ell}_\b{y})$ of the
general position arrangement $\B^{n,\ell}_\b{y}$ may be realized as a
generic $\ell$-dimensional section of the complement
$\sfE^n_\b{y}=M(\B^n_\b{y})$ of the Boolean arrangement $\B^n_\b{y}$,
and the map $\sigma_\b{y}:\sfE^{n,\ll}_\b{y} \to \sfE^n_\b{y}$ as the
inclusion of this generic section.  In the notation established in the
construction of the Morse theoretic complex $K^\bul$ preceding Theorem
\ref{thm:Kdot}, we have $\sfE^{n,\ell}_\b{y} = \FF^\ell \cap
\sfE^n_\b{y}$.  From this construction, it is clear that the complex
$K^\bul(\B^{n,\ell}_\b{y})$ is the rank $\ell$ truncation of the
complex $K^\bul(\B^n_\b{y})$:
\[
K^q(\B^{n,\ell}_\b{y}) =
\begin{cases}
K^q(\B^n_\b{y}) & \text{if $q \le \ell$,}\\
0 & \text{if $q>\ell$,}
\end{cases}
\quad
\text{and}
\quad
\sigma_\b{y}^q =
\begin{cases}
\id:K^q(\B^n_\b{y}) \to K^q(\B^n_\b{y}) & \text{if $q \le \ell$,}\\
0 & \text{if $q>\ell$.}
\end{cases}
\]
See \cite[\S7]{C1} and \cite[Ex.~2.7]{CO1} for explicit constructions
of these complexes.
\end{rem}

The following is a consequence of the above discussion.

\begin{prop} \label{prop:gpPhi}
The representations $\Phi^\bul_\cG:\pi_1(\sfB(\cG),\b{y}) \to
\Aut_\C(K^\bul(\B^{n,\ell}_\b{y}))$ and
$\widetilde\Phi^\bul_\cG:\pi_1(\sfB(\cG),\b{y}) \to
\Aut_\C(K^\bul(\B^{n}_\b{y}))$ satisfy $\Phi^q_\cG =
\widetilde\Phi^q_\cG$ for $q \le \ell$.
\end{prop}

Now recall the bundle $p_\cT:\sfE(\cT) \to \sfY(\cT)$, with fiber
$\sfE^n_\b{x}=M(\B^n_\b{x})$, and the inclusion map $j_\cT:\sfB(\cT)
\to \sfY(\cT)$ from Section \ref{sec:moduli}.  As discussed above, an
element $\c \in \pi_1(\sfB(\cT),\b{x})$ induces a chain automorphism
$\widetilde\Phi^\bul_\cT(\c)$ of the cochain complex
$K^\bul(\B^n_\b{x})$.  Let $g:\bS^1 \to \sfB(\cT)$ be a loop
representing $\c$.  Then $j_\cT \circ g:\bS^1 \to \sfY(\cT)$ is a loop
representing an element, say $\xi \in \pi_1(\sfY(\cT),\b{x})$, and
$\xi$ also induces a chain automorphism of the complex
$K^\bul(\B^n_\b{x})$, which we denote by $\Xi^\bul(\c)$.

\begin{lem}\label{lem:same}
The chain automorphisms $\widetilde\Phi^\bul_\cT(\c)$ and
$\Xi^\bul(\c)$ of the complex $K^\bul(\B^n_\b{x})$ coincide,
$\widetilde\Phi^\bul_\cT(\c) = \Xi^\bul(\c)$.
\end{lem}
\begin{proof}
Let $g^*(p'_\cT)$ denote the pullback of the bundle $p'_\cT:\sfE'(\cT)
\to \sfB(\cT)$ along the map $g:\bS^1 \to \sfB(\cT)$, and let
$(j_\cT\circ g)^*(p_\cT)$ be the pullback of $p_\cT:\sfE(\cT) \to
\sfY(\cT)$ along $j_\cT\circ g:\bS^1 \to\sfY(\cT)$.  The chain
automorphisms $\widetilde\Phi^\bul_\cT(\c)$ and $\Xi^\bul(\c)$ are
induced by the monodromies of the bundles $g^*(p'_\cT)$ and
$(j_\cT\circ g)^*(p_\cT)$ respectively.  Checking that these two
bundles are identical, we conclude that the two monodromies are the
same.  The result follows.
\end{proof}

For any combinatorial type $\cT$, there is an inclusion
$i_\cT:\sfB(\cG) \to \sfY(\cT)$ of the moduli space
$\sfB(\cG)=\sfY(\cG)$ of all general position arrangements in the
space $\sfY(\cT)$.  More generally, if $\cS$ and $\cT$ are types
satisfying $\ind(\cT) \subseteq \ind(\cS)$, there is an inclusion
$i_{\cS,\cT}:\sfY(\cS) \to \sfY(\cT)$.

\begin{prop}\label{prop:onto}
If $\cS$ and $\cT$ are combinatorial types with $\ind(\cT) \subseteq
\ind(\cS)$, the map $(i_{\cS,\cT})_*:H_1(\sfY(\cS)) \to
H_1(\sfY(\cT))$ induced by the inclusion $i_{\cS,\cT}:\sfY(\cS) \to
\sfY(\cT)$ is surjective.
\end{prop}
\begin{proof}
It suffices to consider the case where $\ind(\cS)=\ind(\cT) \cup \{J\}$
for some $J=\{j_1,\dots,j_{\ll+1}\}\subset [n+1]$ with $J \notin
\ind(\cT)$.  Let $\cT_J$ be the type with $\ind(\cT_J)=\{J\}$, and let
$U=\sfY(\cT)$ and $V=\sfY(\cT_J)=(\CP^\ll)^n \setminus
\{\Delta_J=0\}$.  Then $U \cap V = \sfY(\cS)$, and $U \cup V =
(\CP^\ll)^n \setminus \bigcup_{q=1}^k \left(\{\Delta_{I_q}=0\} \cap
\{\Delta_J=0\}\right)$, where $\ind(\cT)=\{I_1,\dots,I_k\}$.  Since
$\{\Delta_{I_q}=0\} \cap \{\Delta_J=0\}$ is of complex codimension two
in $(\CP^\ll)^n$ for each $q$, the inclusion of $U\cup V$ in
$(\CP^\ll)^n$ induces an isomorphism $H_1(U\cup V) \xrightarrow{\sim}
H_1((\CP^\ll)^n)$.  Thus $H_1(U \cup V)=0$.  With these observations,
the Mayer-Vietoris sequence in homology is of the form
\[
\dots \to H_2(U\cup V) \to H_1(\sfY(\cS)) \to
H_1(\sfY(\cT)) \oplus H_1(\sfY(\cT_J)) \to H_1(U \cup V)= 0.
\]
Since the map $H_1(\sfY(\cS)) \twoheadrightarrow H_1(\sfY(\cT)) \oplus
H_1(\sfY(\cT_J))$ is given by
$\left((i_{\cS,\cT})_*,(i_{\cS,\cT_J})_*\right)$, the result follows.
\end{proof}

\begin{cor} \label{cor:GPonto}
For any combinatorial type $\cT$, the inclusion $i_\cT:B(\cG) \to
\sfY(\cT)$ induces a surjection $(i_\cT)_*:H_1(\sfB(\cG)) \to
H_1(\sfY(\cT))$.
\end{cor}

The above results provide a means for comparing the representations
\[
\Phi^\bul_\cT:\pi_1(\sfB(\cT),\b{x}) \to \Aut_\C(K^\bul(\A_\b{x}))
\quad \text{and} \quad \Phi^\bul_\cG:\pi_1(\sfB(\cG),\b{y})\to
\Aut_C(K^\bul(\B^{n,\ll}_\b{y}),
\]
corresponding to an arbitrary arrangement of $n$ hyperplanes in
$\C^\ll$ and a general position arrangement of $n$ hyperplanes in
$\C^\ll$.

\begin{thm} \label{thm:conj}
Let $\cT$ be a combinatorial type with $I_0 \in \ind(\cT)$.  Let
$\A_\b{x}$ be an arrangement of type $\cT$, corresponding to a point
$\b{x}\in \sfB(\cT)$.  Let $\B^{n,\ll}_\b{y}$ be a
general position arrangement, corresponding to a point $\b{y}\in
\sfB(\cG)$.  Then for any element $\c \in
\pi_1(\sfB(\cT),\b{x})$, there is an element $\Gamma \in
\pi_1(\sfB(\cG),\b{y})$ so that the diagram
\begin{equation*} \label{eq:conjcd}
\begin{CD}
K^{\bul}(\B^{n,\ll}_{\b{y}}) @>s_{\b{x}}^{\bul}>> K^{\bul}(\A_{\b{x}})\\
@VV\Phi^{\bul}_\cG(\Gamma)V      @VV\Phi^{\bul}_\cT(\c)V \\
K^{\bul}(\B^{n,\ll}_{\b{y}}) @>s_{\b{x}}^{\bul}>> K^{\bul}(\A_{\b{x}})
\end{CD}
\end{equation*}
commutes up to conjugacy.
\end{thm}
\begin{proof}
Let $g:\bS^1 \to \sfB(\cT)$ be a loop representing
$\c\in\pi_1(\sfB(\cT),\b{x})$.  Then $j_\cT \circ g:\bS^1 \to
\sfB(\cT) \to \sfY(\cT)$ is a loop representing $\xi \in
\pi_1(\sfY(\cT),\b{x})$.  Using Proposition \ref{prop:Acd} and Lemma
\ref{lem:same}, we have the commutative diagram
\[
\begin{CD}
K^{\bul}(\B^{n}_{\b{x}}) @>\id>> K^{\bul}(\B^{n}_{\b{x}})
@>s_{\b{x}}^{\bul}>> K^{\bul}(\A_{\b{x}})\\
@VV{\Xi^\bul(\c)}V      @VV\widetilde\Phi^{\bul}_\cT(\c)V
@VV\Phi^{\bul}_\cT(\c)V \\
K^{\bul}(\B^{n}_{\b{x}}) @>\id>> K^{\bul}(\B^{n}_{\b{x}})
@>s_{\b{x}}^{\bul}>> K^{\bul}(\A_{\b{x}})
\end{CD}
\]
where $\Xi^\bul(\c)$ is the chain automorphism induced by $\xi$.

Denote the homology class of $\xi$ in $H_1(\sfY(\cT))$ by $[\xi]$.  By
Corollary \ref{cor:GPonto}, the inclusion
$i_\cT:\sfB(\cG)\to\sfY(\cT)$ induces a surjection
$(i_\cT)_*:H_1(\sfB(\cG)) \to H_1(\sfY(\cT))$.  Let
$\b{y}\in\sfB(\cG)$ be a basepoint (the point $\b{x}\in\sfB(\cT)$ is
not in $\sfB(\cG)$ if $\cT \neq \cG$).  Let $\Gamma \in
\pi_1(\sfB(\cG),\b{y})$ be an element whose homology class $[\Gamma]$
satisfies $(i_\cT)_*([\Gamma]) = [\xi]$.  Then the chain automorphisms
induced by $(i_\cT)_\#(\Gamma) \in \pi_1(\sfY(\cT),\b{y})$ and $\xi
\in \pi_1(\sfY(\cT),\b{x})$ are conjugate, so we have the commutative
diagram
\[
\begin{CD}
K^{\bul}(\B^{n}_{\b{y}}) @>\sim>> K^{\bul}(\B^{n}_{\b{x}})
@>\id>> K^{\bul}(\B^{n}_{\b{x}})
@>s_{\b{x}}^{\bul}>> K^{\bul}(\A_{\b{x}})\\
@VV{(i_\cT)_\#(\Gamma)^\bul}V @VV{\Xi^\bul(\c)}V
@VV\widetilde\Phi^{\bul}_\cT(\c)V   @VV\Phi^{\bul}_\cT(\c)V \\
K^{\bul}(\B^{n}_{\b{y}}) @>\sim>> K^{\bul}(\B^{n}_{\b{x}})
@>\id>> K^{\bul}(\B^{n}_{\b{x}})
@>s_{\b{x}}^{\bul}>> K^{\bul}(\A_{\b{x}})
\end{CD}
\]

Recall that $\sfB(\cG)=\sfY(\cG)$.  Over $\sfB(\cG)$ and $\sfY(\cT)$,
we have bundles $p_\cG:\sfE(\cG) \to\sfB(\cG)$ and $p_\cT:\sfE(\cT)
\to\sfY(\cT)$, each with fiber $\sfE^n_\b{y}$, the complement of the
Boolean arrangement $\B^n_\b{y}$.  Let $G:\bS^1 \to \sfB(\cG)$ be a
loop representing $\Gamma \in \pi_1(\sfB(\cG),\b{y})$.  It is then
readily checked that the pullbacks $G^*(p_\cG)$ and $(i_\cT \circ
G)^*(p_\cT)$ of the two bundles above are identical.  Consequently,
the chain automorphisms $\widetilde\Phi^{\bul}_\cG(\Gamma)$ and
$(i_\cT)_\#(\Gamma)^\bul$ of $K^\bul(\B^n_\b{y})$ induced by $\Gamma$
and $(i_\cT)_\#(\Gamma)$ are equal.  This fact, together with
\eqref{eq:Bcd}, yields the commutative diagram
\[
\begin{CD}
K^{\bul}(\B^{n,\ll}_{\b{y}}) @>{\sigma^\bul_\b{y}}>>
K^{\bul}(\B^{n}_{\b{y}}) @>\id>> K^{\bul}(\B^{n}_{\b{y}})
@>\sim>>  K^{\bul}(\B^{n}_{\b{x}})
@>s_{\b{x}}^{\bul}>> K^{\bul}(\A_{\b{x}})\\
@VV\Phi^{\bul}_\cG(\Gamma)V  @VV\widetilde\Phi^{\bul}_\cG(\Gamma)V
@VV{(i_\cT)_\#(\Gamma)^\bul}V  @VV\widetilde\Phi^{\bul}_\cT(\c)V
@VV\Phi^{\bul}_\cT(\c)V \\
K^{\bul}(\B^{n,\ll}_{\b{y}}) @>{\sigma^\bul_\b{y}}>>
K^{\bul}(\B^{n}_{\b{y}}) @>\id>> K^{\bul}(\B^{n}_{\b{y}})
@>\sim>>  K^{\bul}(\B^{n}_{\b{x}}) @>s_{\b{x}}^{\bul}>>
K^{\bul}(\A_{\b{x}})
\end{CD}
\]
Since $\sigma^q_\b{y}=\id:K^q(\B^{n,\ll}_\b{y}) \xrightarrow{=}
K^q(\B^n,\b{y})$ and
$\Phi^{\bul}_\cG(\Gamma)=\widetilde\Phi^{\bul}_\cG(\Gamma)$ for $q \le
\ll$ by Remark~\ref{rem:gp} and Proposition \ref{prop:gpPhi}, and
$K^q(\B^{n,\ll}_\b{y})=K^q(\A_\b{x})=0$ for $q>\ll$, this completes
the proof.
\end{proof}

Given $\c\in\pi_1(\sfB(\cT),\b{x})$, we say that $\Gamma \in
\pi_1(\sfB(\cG),\b{y})$ is {\em compatible} with $\c$ if their
homology classes satisfy $(i_\cT)_*([\Gamma])=(j_\cT)_*([\c])$.  The
results established above show that for every $\c$ there exists a
compatible $\Gamma$.  A key tool in this paper is the relationship
between the automorphism $\Phi^\bul_\cT(\c) \in \Aut_C
K^\bul(\A_\b{x})$ for an {\em arbitrary} arrangement and the
automorphism $\Phi^\bul_\cG(\Gamma) \in \Aut_\C
K^\bul(\B^{n,\ll}_\b{y})$ for a {\em general position} arrangement
provided by Theorem \ref{thm:conj} for compatible $\c$ and $\Gamma$.
This relationship extends to cohomology.  By Theorem
\ref{thm:inducedrep}, the representation $\Psi^q_\cT:
\pi_1(\sfB(\cT),\b{x}) \to \Aut_\C(H^q(\sfM_\b{x};\LL_\b{x}))$ in
cohomology is induced by the representation
$\Phi^\bul_\cT:\pi_1(\sfB(\cT),\b{x}) \to \Aut_\C(K^\bul(\A_\b{x}))$.
Let $\bar{s}^q_\b{x}:H^q(\sfE^{n.\ll}_\b{y};\LL_\b{y}) \to
H^q(\sfM_\b{x};\LL_\b{x})$ denote the map in cohomology induced by
$s_{\b{x}}^{\bul}:K^{\bul}(\B^{n,\ll}_{\b{y}}) \to
K^{\bul}(\A_{\b{x}})$.  If $\Gamma\in\pi_1(\sfB(\cG),\b{y})$ is
compatible with $\c \in\pi_1(\sfB(\cT),\b{x})$, then Theorem
\ref{thm:conj} implies that, up to conjugacy, we have $\bar{s}^q_\b{x}
\circ \Psi^q_\cG(\Gamma) = \Psi^q_\cT(\c) \circ \bar{s}^q_\b{x}$ for
each $q$, $0\le q \le \ll$.  We will pursue the implications of this
relationship for
Gauss-Manin connections in Section \ref{sec:GM}.

\section{Compatible Classes}
\label{sec:compatible}
For the type $\cG$ of general position arrangements, the closure of
the moduli space is $\overline{\sfB}(\cG)= (\CP^\ll)^n$.  Recall that
points in $(\CP^\ll)^n$ are given by matrices $\b{x}$ of the form
\eqref{eq:point}, and that $\D_J(\b{x})$ denotes the determinant of
the submatrix of $\b{x}$ with rows specified by $\ll +1$ element
subsets $J=\{j_1,\dots,j_{\ll+1}\}$ of $[n+1]=\{1,\dots,n+1\}$.  The
divisor $\sfD(\cG)=\overline{\sfB}(\cG)\setminus \sfB(\cG)$ is given
by $\sfD(\cG) = \bigcup_J \sfD_J$, whose components,
$\sfD_J =\{\b{x} \in (\CP^\ll)^n \mid \D_J(\b{x})=0\}$, are
irreducible hypersurfaces indexed by $J$.

Choose a basepoint $\b{y} \in \sfB(\cG)$, and for each $\ll+1$ element
subset $J$ of $[n+1]$, let $\b{z}_J$ be a generic point in $\sfD_J$.
Let $\Gamma_{\!J}$ be a meridian loop based at $\b{y}$ in $\sfB(\cG)$
about the point $\b{z}_J \in \sfD_J$.  Note that $\b{y} \in \sfY(\cT)$
and that $\Gamma_{\!J}$ is a (possibly null-homotopic) loop in
$\sfY(\cT)$ for any combinatorial type $\cT$.

\begin{prop} \label{prop:H1YT}
For any combinatorial type $\cT$, the homology group $H_1(\sfY(\cT))$
is generated by the classes $\{[\Gamma_{\!J}] \mid J \in \ind(\cT)\}$.
\end{prop}
\begin{proof}
Given $\cT$, let $Z = \bigcup_{J \in \ind(\cT)} \sfD_J$,
so that $\sfY(\cT) = (\CP^\ll)^n \setminus Z$.  Denote by $\Sigma{Z}$
the singular set of $Z$, and let $X=(\CP^\ll)^n \setminus \Sigma{Z}$
and $D=Z \setminus\Sigma{Z}$.  Consider the corresponding Gysin
sequence in homology with integer coefficients, as discussed in
\cite[p.~46]{Di}.  This sequence is of the form
\[
\dots \to H_{k-1}(D) \xrightarrow{\theta} H_k(X \setminus D)
\xrightarrow{i_*} H_k(X) \to H_{k-2}(D) \xrightarrow{\theta}
H_{k-1}(X \setminus D) \to \dots
\]
where $i:X \setminus D \to X$ denotes the inclusion.

Since the complex codimension of $\Sigma{Z}$ is at least two in
$(\CP^\ll)^n$, the inclusion of $X$ in $(\CP^\ll)^n$ induces an
isomorphism $H_1(X) \simeq H_1((\CP^\ll)^n$ and thus $H_1(X)=0$.  So
for small $k$, the above sequence is of the form \[ \dots \to H_2(X)
\to H_0(D) \xrightarrow{\theta} H_1(X \setminus D) \to H_1(X)=0 \] The
connected components of $D=Z \setminus\Sigma{Z}$ are in one-to-one
correspondence with the irreducible components $\sfD_J$ of $Z$.
Hence, $H_0(D)$ is free abelian, of rank equal to the cardinality of
$\ind(\cT)$.  A basis for $H_0(D)$ is given by the classes
$[\b{z}_J]$, $J\in \ind(\cT)$, where $\b{z}_J \in \sfD_J$ is the
generic point in $\sfD_J$ chosen above.  Thus the classes
$\theta([\b{z}_J])$ generate $H_1(X \setminus D) = H_1(\sfY(\cT))$.

Finally, let $T$ be a tubular neighborhood of $D$, and let
$p:\partial{T} \to D$ denote the projection from the boundary of $T$
to $D$.  Then the map $\theta:H_0(D) \to H_1(X \setminus D)$ is given
by $\theta([\b{z}])=[p^{-1}(\b{z})]$.  Since $[p^{-1}(\b{z}_J)]$
clearly coincides with $[\Gamma_{\!J}]$ up to sign for each $J$, this
completes the proof.
\end{proof}

\begin{cor} \label{cor:H1BG}
The homology group $H_1(\sfB(\cG))$ is generated by the classes
$[\Gamma_{\!J}]$, where $J$ ranges over all $\ll+1$ element subsets of
$[n+1]$.
\end{cor}

For a combinatorial type $\cT$, recall that $\dep(\cT)$ denotes the
set of all $\ll+1$ element subsets $J=\{j_1,\dots,j_{\ll+1}\}$ of
$[n+1]$ for which the determinant $\Delta_J(\b{x})$ vanishes for
$\b{x} \in \sfB(\cT)$.  If $\cT$ is realizable, $\dep(\cT)$ is the set
of all $J$ for which $\{H_{j_1},\dots,H_{j_{\ll+1}}\}$ is linearly
dependent in the projective closure of an arrangement $\A$ of type
$\cT$.  Impose a partial order on combinatorial types as follows: $\cT
\ge \cT' \iff \dep(\cT) \subseteq \dep(\cT')$.  Note that the
combinatorial type $\cG$ of general position arrangements is the
maximal element with respect to this partial order.

Write $\cT > \cT'$ if $\dep(\cT) \subsetneq \dep(\cT')$.  In this case
we define the relative dependence set $\dep(\cT',\cT)= \dep(\cT')
\setminus \dep(\cT)$.  If $\cT > \cT'$, we say that $\cT$ {\em covers}
$\cT'$ and $\cT'$ is a {\em degeneration} of  $\cT$
if there is no combinatorial type $\cT''$ with
$\cT>\cT''>\cT'$.

\begin{lem} \label{lem:cover}
The moduli space $\sfB(\cT')$ has complex codimension one in the
closure $\overline{\sfB}(\cT)$ of the moduli space $\sfB(\cT)$ if and
only if $\cT$ covers $\cT'$.
\end{lem}
\begin{proof}
If $\sfB(\cT')$ is codimension one in $\overline{\sfB}(\cT)$, then
clearly $\cT$ covers $\cT'$.  Conversely, if $\cT$ covers $\cT'$, then
$\sfB(\cT') \subset \overline{\sfB}(\cT)$ is defined by the vanishing
of $\Delta_J$ for any $J \in \dep(\cT',\cT)$.
\end{proof}

Recall that $j_\cT:\sfB(\cT) \to \sfY(\cT)$ denotes the natural
inclusion.

\begin{thm} \label{thm:lincomb}
Let $\cT$ be a combinatorial type which covers the type $\cT'$.  Let
$\b{x}'$ be a point in $\sfB(\cT')$, and $\c \in
\pi_1(\sfB(\cT),\b{x})$ a simple loop in $\sfB(\cT)$ about $\b{x}'$.
Then the homology class $[\c]$ satisfies
\begin{equation} \label{eq:compat}
(j_\cT)_*([\c]) = \sum_{J \in \dep(\cT',\cT)} m_J \cdot
[\Gamma_{\!J}],
\end{equation}
where $m_J$ is the order of vanishing of the restriction of $\Delta_J$
to $\overline{\sfB}(\cT)$ along $\sfB(\cT')$.
\end{thm}
\begin{proof}
By Proposition \ref{prop:H1YT}, the classes
$\{[\Gamma_{\!J}] \mid J \in \ind(\cT)\}$ generate $H_1(\sfY(\cT))$.
So the image of $[\c] \in H_1(\sfB(\cT))$ in $H_1(\sfY(\cT))$ is of
the form $(j_\cT)_*([\c]) = \sum m_J \cdot [\Gamma_{\!J}]$, where the
sum is over all $J \in \ind(\cT)$, and $m_J \in \Z$.

To determine the coefficients $m_J$, let $f=0$ be a local defining
equation for $\sfB(\cT')$ in $\overline{\sfB}(\cT)$ near $\b{x}$.  By
hypothesis, the winding number of $f$ about $\c$ is one:
\[
\frac{1}{2 \pi \ii} \int_\c \frac{df}{f} = 1.
\]
For each $J \in \ind(\cT)$, let $B_J$ be a
disk with $\partial B_J = \c$, $B_J \cap \sfB(\cT') = \b{x}'$, and
$B_J \setminus \{\b{x}'\} \subset \sfB(\cT)$.  Suppose also that
$F_J:(\CP^\ll)^n \to\C$ defines the divisor $\sfD_J=\{\Delta_J=0\}$
locally near $\b{x}$, and let $\overline{j}_\cT:\overline{\sfB}(\cT)
\to (\CP^\ll)^n$ denote the inclusion of the closure of $\sfB(\cT)$ in
$(\CP^\ll)^n$.  The coefficient $m_J$ is then given by the winding
number of $F_J \circ\overline{j}_\cT$ about $\c$:
\[
m_J=\frac{1}{2 \pi \ii} \int_\c \frac{d(F_J
\circ\overline{j}_\cT)}{F_J \circ\overline{j}_\cT}.
\]

If $J\in \ind(\cT') \cap \ind(\cT)$, then $F_J(\b{x}') \neq 0$.  So
for such $J$, the restriction of $F_J$ to the disk $B_J$ is never
zero.  Consequently, we have $\int_\c {d(F_J
\circ\overline{j}_\cT)}/{F_J \circ\overline{j}_\cT} = 0$, and thus
$m_J = 0$ if $J \in \ind(\cT')$.

If $J \in \dep(\cT',\cT)$, then the restriction, $F_J
\circ\overline{j}_\cT$, of $F_J$ to $\overline{\sfB}(\cT)$ defines
a hypersurface which contains $\sfB(\cT')$.  So near $\b{x}$, we have
$F_J \circ \overline{j}_\cT = c f^{m}$, where $c \in \C^*$, and $m \in
\Z_{>0}$ is the order of vanishing of $F_J \circ\overline{j}_\cT$
along $\sfB(\cT')$.  Hence for $J \in \dep(\cT',\cT)$,
the coefficient
\[
m_J=\frac{1}{2 \pi \ii} \int_\c \frac{d(F_J
\circ\overline{j}_\cT)}{F_J \circ\overline{j}_\cT} =\frac{1}{2 \pi
\ii} \int_\c \frac{d(c f^m)}{c f^m} =\frac{m}{2 \pi \ii} \int_\c
\frac{df}{f} =m
\]
is as asserted.
\end{proof}

Recall that we say $\Gamma \in \pi_1(\sfB(\cG),\b{y})$ is compatible
with $\c\in\pi_1(\sfB(\cT),\b{x})$ if the homology classes $[\Gamma]$
and $[\c]$ satisfy $(i_\cT)_*([\Gamma])=(j_\cT)_*([\c])$.

\begin{cor} \label{cor:compatible}
If $\cT$ covers $\cT'$ and $\c \in \pi_1(\sfB(\cT),\b{x})$ is a simple
loop in $\sfB(\cT)$ about $\b{x}' \in\sfB(\cT')$, then any loop
$\Gamma \in \pi_1(\sfB(\cG),\b{y})$ whose homology class satisfies the
condition $[\Gamma]=\sum_J m_J \cdot [\Gamma_{\!J}]$ given by
\eqref{eq:compat} is compatible with $\c$.
\end{cor}

\section{Nonresonant Weights} \label{sec:weights}
If $\LL$ is a nontrivial rank one local system on the complement
$\sfE^{n,\ll}$ of a general position arrangement $\B^{n,\ll}$ in
$\C^\ll$, then the cohomology $H^q(\sfE^{n,\ll};\LL)$ vanishes in all
dimensions $q$ except possibly $q=\ll$.  See \cite[\S7]{C1} for a
proof of this fact using the complex $K^\bul(\B^{n,\ll})$.  For an
arbitrary arrangement $\A$, a collection of weights $\la$ and the
corresponding local system $\LL$ are called {\em nonresonant} if the
cohomology $H^q(\sfM;\LL)$ vanishes in all dimensions $q$ except
possibly $q=\ll$.  We now recall a combinatorial condition due to
Schechtman, Terao, and Varchenko \cite{STV} which insures
nonresonance.  An edge of an arrangement $\A$ is a nonempty
intersection of hyperplanes in $\A$.  An edge is called {\em dense} if
the subarrangement of hyperplanes containing it is irreducible: the
hyperplanes cannot be partitioned into nonempty sets so that, after a
change of coordinates, hyperplanes in different sets are in different
coordinates.  This is a combinatorially determined property, see
\cite{STV}.  For each edge $X$, define $\la_X=\sum_{X \subseteq
H_j}\la_j$.  Let $\Ai=\A\cup H_{n+1}$ be the projective closure of
$\A$, the union of $\A$ and the hyperplane at infinity in $\CP^\ll$,
see \cite{OT2}.  Set $\la_{n+1} = -\sum_{j=1}^n \la_j$.

\begin{thm}[\cite{STV}] \label{thm:vanish}
Let $\sfM$ be the complement of an essential arrangement $\A$ in
$\C^\ll$ of combinatorial type $\cT$.  If $\LL$ is a rank one local
system on $\sfM$ whose weights $\la$ satisfy the condition\\[4pt]
{\rm (STV) \hfill $\la_X \notin \Z_{\ge 0}$ for every dense edge
$X$ of $\Ai$,\hfill \phantom{{\rm (STV)}}}\\[4pt]
then $H^q(\sfM;\LL) = 0$ for $q\neq\ll$ and $\dim H^\ll(\sfM;\LL) =
|\chi(\sfM)|$, where $\chi(\sfM)$ is the Euler characteristic of
$\sfM$.
\end{thm}
Thus the (STV)-condition implies nonresonance.  In this paper, we call
local systems $\LL$ whose weights $\la$ satisfy the (STV)-condition
{\em $\cT$--nonresonant}.

\begin{rem} \label{rem:gpnr}
The hyperplanes are among the dense edges for any arrangement, and are
the only dense edges for a general position arrangement.
Consequently, if $\la$ is $\cT$--nonresonant, then $\la$ is
$\cG$--nonresonant.
\end{rem}

Any collection of weights $\la$ determines an element $a_\la
=\sum_{j=1}^n \la_j\, a_j\in A^1(\A)$ in the Orlik-Solomon algebra of
$\A$.  Since $a_\la \wedge a_\la=0$, multiplication by $a_\la$ defines
a differential on $A(\A)$.  The resulting complex
$(A^\bul(\A),a_\la\wedge)$ may be identified with a subcomplex of the
twisted de Rham complex of $\sfM$ with coefficients in $\LL$.  Theorem
\ref{thm:vanish} above may be established by showing that, for
$\cT$--nonresonant weights $\la$, there is an isomorphism
$H^*(\sfM;\LL) \simeq H^*(A^\bul(\A),a_\la\wedge)$ \cite{ESV,STV}, and
that $H^q(A^\bul(\A),a_\la\wedge)=0$ for $q\neq\ll$ and $\dim
H^\ll(A^\bul(\A),a_\la\wedge)=|\chi(\sfM)|$, see \cite{Yuz}.

For $\cT$--nonresonant weights, one can exhibit an explicit basis, the
$\beta${\bf{nbc}} basis of \cite{FT}, for the single nonvanishing
cohomology group $H^\ll(\sfM;\LL) \simeq
H^\ll(A^\bul(\A),a_\la\wedge)$.  Recall that the hyperplanes of
$\A=\{H_j\}_{j=1}^n$ are ordered.  A circuit is an inclusion-minimal
dependent set of hyperplanes in $\A$, and a broken circuit is a set
$S$ for which there exists $H < \min(S)$ so that $S \cup\{H\}$ is a
circuit.  A frame is a maximal independent set, and an {\bf nbc} frame
is a frame which contains no broken circuit.  Since $\A\subset\C^\ll$
is essential, every frame has cardinality~$\ll$.  An {\bf{nbc}} frame
$B=(H_{j_1},\dots,H_{j_\ll})$ is a $\beta${\bf{nbc}} frame provided
that for each $k$, $1\le k \le \ll$, there exists $H\in \A$ such that
$H<H_{j_k}$ and $(B\setminus\{H_{j_k}\})\cup\{H\}$ is a frame.

Let $\beta${\bf{nbc}}$(\A)$ be the set of all $\beta${\bf{nbc}} frames
of $\A$.  Given $B=(H_{j_1},\dots,H_{j_\ll})$ in
$\beta\text{\bf{nbc}}(\A)$, define $\zeta(B) \in A^\ll(\A)$ by
$\zeta(B) = \wedge_{p=1}^\ll a_\la(X_p)$, where $X_p=\bigcap_{k=p}^\ll
H_{j_k}$ and $a_\la(X) = \sum_{X \subseteq H_i} \la_{i} a_{i}$.
Denote the cohomology class of $\zeta(B)$ in
$H^\ll(A^\bul(\A),a_\la\wedge)$ by the same symbol.

\begin{thm}[\cite{FT}]\label{thm:bnbc}
Let $\sfM$ be the complement of an essential arrangement $\A$ in
$\C^\ll$ of combinatorial type $\cT$.  If $\LL$ is a rank one local
system on $\sfM$ corresponding to $\cT$--nonresonant weights $\la$,
then the set $\{\zeta(B) \mid B \in \beta\text{\bf{nbc}}(\A)\}$ is a
basis for the only nontrivial local system cohomology group
$H^\ll(\sfM;\LL) \simeq H^\ll(A^\bul(\A),a_\la\wedge)$.
\end{thm}

\begin{exm} \label{exm:GPbnbc}
If $\B^{n,\ll}$ is a general position arrangement of $n$ hyperplanes
in $\C^\ll$, then
$\beta\text{\bf{nbc}}(\B^{n,\ll})=\{(H_{j_1},\dots,H_{j_\ll}) \mid 2
\le j_1 < \dots < j_\ll \le n\}$.  For $\cG$--nonresonant weights
$\la$, it follows that the $\beta${\bf{nbc}} basis for
$H^\ll(\sfE^{n,\ll};\LL)$ consists of monomials $\eta_T =
\la_{j_1}\cdots \la_{j_\ll} a_T$, where $T=\{j_1,\dots,j_\ll\}$, $2
\le j_1 < \dots < j_\ll \le n$, and $a_T=a_{j_1}\wedge \dots \wedge
a_{j_\ll}$.
\end{exm}

\begin{thm} \label{thm:nonresprojection}
For $\cT$--nonresonant weights $\la$, the map $s:\sfM \to \sfE^n$
induces a surjection $s^*:H^\ll(\sfE^{n,\ll},\LL) \twoheadrightarrow
H^\ll(\sfM;\LL)$.
\end{thm}
\begin{proof}
The map $s:\sfM \to \sfE^n$ induces a surjection $s^*:A(\B^{n,\ll})
\twoheadrightarrow A(\A)$ by Proposition \ref{prop:GPsurjection}.
This is a chain map $s^*:(A^\bul(\B^{n,\ll}),a_\la\wedge) \to
(A^\bul(\A),a_\la\wedge)$, as is readily checked.  So there is an
induced map $\bar{s}:H^\ll(A^\bul(\B^{n,\ll}),a_\la\wedge) \to
H^\ll(A^\bul(\A),a_\la\wedge)$ in cohomology.

By Remark \ref{rem:gpnr}, if $\la$ is $\cT$--nonresonant, then $\la$
is $\cG$--nonresonant.  For such weights, we have isomorphisms
$H^\ll(\sfE^{n,\ll};\LL) \simeq H^\ll(A^\bul(\B^{n,\ll}),a_\la\wedge)$
and $H^\ll(\sfM;\LL) \simeq H^\ll(A^\bul(\A),a_\la\wedge)$.
Furthermore, the elements of the $\beta${\bf{nbc}} basis for
$H^\ll(A^\bul(\A),a_\la\wedge)$ are linear combinations of the
monomials $\eta_T$, which form the $\beta${\bf{nbc}} basis for
$H^\ll(A^\bul(\B^{n,\ll}),a_\la\wedge)$ as noted in Example
\ref{exm:GPbnbc}.  This fact, together with the surjectivity of the
map $s^*:A(\B^{n,\ll}) \twoheadrightarrow A(\A)$ completes the proof.
\end{proof}
This leads to the following consequence of Theorem \ref{thm:conj}.

\begin{cor} \label{cor:Hconj}
Let $\la$ be a collection of $\cT$--nonresonant weights.  If
$\c\in\pi_1(\sfB(\cT),\b{x})$ and $\Gamma \in \pi_1(\sfB(\cG),\b{y})$
are compatible, then the conjugacy class of the automorphism
$\Psi^\ll_\cT(\c) \in \Aut_\C(H^\ll(\sfM_\b{x};\LL_\b{x}))$ is
determined by $\Psi^\ll_\cG(\Gamma) \in
\Aut_\C(H^\ll(\sfE^{n,\ll}_\b{y};\LL_\b{y}))$ and the surjection
$\bar{s}^\ll_\b{x}:H^\ll(\sfE^{n,\ll}_\b{y};\LL_\b{y})
\twoheadrightarrow H^\ll(\sfM_\b{x};\LL_\b{x})$.
\end{cor}

\section{Gauss-Manin Connections}
\label{sec:GM}
The vector bundle $\pi^q:\b{H}^q\to\sfB(\cT)$, with fiber
$(\pi^q)^{-1}(\b{x})=H^q(\sfM_\b{x};\LL_\b{x})$ at
$\b{x}\in\sfB(\cT)$, supports a Gauss-Manin connection corresponding
to the representation $\Psi^q_\cT:\pi_1(\sfB(\cT),\b{x}) \to
\Aut_\C(H^q(\sfM_\b{x};\LL_\b{x}))$.  We now study this connection
using results of the previous section.  In light of Corollary
\ref{cor:Hconj}, for an arrangement $\A$ of combinatorial type $\cT$,
we focus on $\cT$--nonresonant weights and the case $q=\ll$.

Over a manifold, there is a well known equivalence between local
systems and complex vector bundles equipped with flat connections, see
\cite{De,Ko}.  Let $\b{V}\to X$ be such a bundle, with connection
$\nabla$.  The latter is a $\C$-linear map $\nabla:\cE^0(\b{V}) \to
\cE^1(\b{V})$, where $\cE^p(\b{V})$ denotes the complex $p$-forms on
$X$ with values in $\b{V}$, which satisfies $\nabla(f\sigma)= \sigma
df + f \nabla(\sigma)$ for a function $f$ and $\sigma\in\cE^0(\b{V})$.
The connection extends to a map $\nabla:\cE^p(\b{V}) \to
\cE^{p+1}(\b{V})$ for $p\ge 0$, and is flat if the curvature
$\nabla\circ\nabla$ vanishes.  Call two connections $\nabla$ and
$\nabla'$ on $\b{V}$ isomorphic if $\nabla'$ is obtained from $\nabla$
by a gauge transformation, $\nabla'=g\circ\nabla\circ g^{-1}$ for some
$g:X\to\Hom(\b{V},\b{V})$.

The aforementioned equivalence is given by $(\b{V},\nabla) \mapsto
\b{V}^{\nabla}$, where $\b{V}^{\nabla}$ is the local system, or
locally constant sheaf, of horizontal sections $\{\sigma \in
\cE^0(\b{V})\mid \nabla(\sigma)=0\}$.  There is also a well known
equivalence between local systems on $X$ and finite dimensional
complex representations of the fundamental group of $X$.  Note that
isomorphic connections give rise to the same representation.  Under
these equivalences, the local system on $X=\sfB(\cT)$ induced by the
representation $\Psi^q_\cT$ corresponds to a flat connection on the
vector bundle $\pi^q:\b{H}^q\to\sfB(\cT)$, the Gauss-Manin connection.

For $\cT$--nonresonant weights, the $\beta${\bf{nbc}} basis of Theorem
\ref{thm:bnbc} provides a basis for each fiber
$H^\ll(\sfM_\b{x};\LL_\b{x})$, independent of $\b{x}$.  Thus the
vector bundle $\pi^\ll:\b{H}^\ll\to\sfB(\cT)$ is trivial, see
\cite{FT,T1}.  Let $\c\in\pi_1(\sfB(\cT),\b{x})$, and let $g:\bS^1 \to
\sfB(\cT)$ be a representative loop.  Pulling back the trivial bundle
$\pi^\ll:\b{H}^\ll \to \sfB(\cT)$ and the Gauss-Manin connection
$\nabla$, we obtain a flat connection $g^*(\nabla)$ on the trivial
vector bundle over the circle corresponding to the representation of
$\pi_1(\bS^1,1)=\langle\zeta\rangle=\Z$ given by $\zeta \mapsto
\Psi^\ll_\cT(\c)$.  Specifying the flat connection $g^*(\nabla)$
amounts to choosing a logarithm of $\Psi^\ll_\cT(\c)$.  The connection
$g^*(\nabla)$ is determined by a connection $1$-form $dz/z \otimes
\Omega_\cT(\c)$, where the connection matrix $\Omega_\cT(\c)$
corresponding to $\c$ satisfies $\Psi^\ll_\cT(\c) = \exp(-2 \pi\ii
\Omega_\cT(\c))$.  If $\c$ and $\hat\c$ are conjugate in
$\pi_1(\sfB(\cT),\b{x})$, then the resulting connection matrices are
conjugate, and the corresponding connections on the trivial vector
bundle over the circle are isomorphic.  In this sense, the connection
matrix $\Omega_\cT(\c)$ is determined by the homology class $[\c]$ of
$\c$.

\begin{thm} \label{thm:GMCD}
Let $\la$ be a collection of $\cT$--nonresonant weights.  If $\c\in
\pi_1(\sfB(\cT),\b{x})$ and $\Gamma \in \pi_1(\sfB(\cG),\b{y})$ are
compatible, then there is a commutative diagram
\[
\begin{CD}
H^\ll(\sfE^{n,\ll}_\b{y};\LL_\b{y}) @>{\bar{s}^\ll_\b{x}}>>
H^\ll(\sfM_\b{x};\LL_\b{x}) \\
@VV{\Omega_\cG(\Gamma)}V  @VV{\Omega_\cT(\c)}V  \\
H^\ll(\sfE^{n,\ll}_\b{y};\LL_\b{y}) @>{\bar{s}^\ll_\b{x}}>>
H^\ll(\sfM_\b{x};\LL_\b{x})
\end{CD}
\]
Thus a Gauss-Manin connection matrix $\Omega_\cT(\c)$ is determined by
a Gauss-Manin connection matrix $\Omega_\cG(\Gamma)$ and the
surjection $\bar{s}^\ll_\b{x}:H^\ll(\sfE^{n,\ll}_\b{y};\LL_\b{y})
\twoheadrightarrow H^\ll(\sfM_\b{x};\LL_\b{x})$.
\end{thm}
\begin{proof}
Choose a connection matrix $\Omega_\cG(\Gamma)$ so that
$\Psi_\cG^\ll(\Gamma)=\exp(-2\pi\ii \Omega_\cG(\Gamma))$.  Since
$\bar{s}^\ll_\b{x}:H^\ll(\sfE^{n,\ll}_\b{y};\LL_\b{y}) \to
H^\ll(\sfM_\b{x};\LL_\b{x})$ is surjective by Theorem
\ref{thm:nonresprojection}, the endomorphism
$\Omega_\cG(\Gamma):H^\ll(\sfE^{n,\ll}_\b{y};\LL_\b{y}) \to
H^\ll(\sfE^{n,\ll}_\b{y};\LL_\b{y})$ induces an endomorphism
$H^\ll(\sfM_\b{x};\LL_\b{x}) \to H^\ll(\sfM_\b{x};\LL_\b{x})$, which
we denote by $\Omega_\cT(\c)$.  We assert that $\Omega_\cT(\c)$ is a
connection matrix for $\Psi_\cT^\ll(\c)$.  For this, by the above
discussion, it suffices to show that $\exp(-2\pi\ii \Omega_\cT(\c))$
is conjugate to $\Psi_\cT^\ll(\c)$.

By construction, we have $\bar{s}^\ll_\b{x} \circ \Omega_\cG(\Gamma) =
\Omega_\cT(\c) \circ \bar{s}^\ll_\b{x}$.  From this it follows that
$\bar{s}^\ll_\b{x} \circ \exp(-2\pi\ii\Omega_\cG(\Gamma)) =
\exp(-2\pi\ii\Omega_\cT(\c)) \circ \bar{s}^\ll_\b{x}$.  By Theorem
\ref{thm:conj} and Corollary~\ref{cor:Hconj}, up to conjugacy, we have
$\bar{s}^\ll_\b{x} \circ \Psi_\cG^\ll(\Gamma) = \Psi_\cT^\ll(\c) \circ
\bar{s}^\ll_\b{x}$.  Hence, $\exp(-2\pi\ii \Omega_\cT(\c))$ is
conjugate to $\Psi_\cT^\ll(\c)$.
\end{proof}

Next we need formulas of Aomoto and Kita \cite{AK}, see also
\cite{OT2}, for the Gauss-Manin connection matrices in the case of
general position arrangements.  Let $\b{y}\in\sfB(\cG)$, and let
$\B^{n,\ell}_{\b{y}}$ be the corresponding arrangement of type $\cG$,
a general position arrangement of $n$ hyperplanes in $\C^\ell$, with
complement $\sfE^{n,\ell}_{\b{y}}$.  A system of weights
$\la=(\la_1,\dots,\la_n)$ is $\cG$--nonresonant if $\la_j \notin
\Z_{\ge 0}$ for each $j$, $1\le j \le n+1$.  Recall that
$\la_{n+1}=-\sum_{j=1}^n \la_j$.  For such weights, the
$\beta$\textbf{nbc} basis for the cohomology group
$H^\ell(\sfE^{n,\ell}_{\b{y}};\LL_\b{y})$ with coefficients in the
corresponding local system consists of the monomials $\eta_I=
\la_{i_1}\cdots \la_{i_\ll} a_I$, where $I=\{i_1,\dots,i_\ll\}$, $2
\le i_1 < \dots < i_\ll \le n$, and $a_I=a_{i_1}\wedge \dots \wedge
a_{i_\ll}$, see Example \ref{exm:GPbnbc}.

To state the results of Aomoto and Kita on the Gauss-Manin connection
for general position arrangements, we require some notation.  We use
the formulation of \cite[\S{10.3}]{OT2}.  For $J=\{j_1,\dots,j_m\}$,
write $\la_J=\sum_{j\in J}\la_j$, and let $J_p=\{j_1,\dots
\widehat{j_p},\dots,j_m\}$ for each $p$, $1\le p \le m$.  If
$I,I'\subset [n]$, $|I|=|I'|=\ell$, and $|I\cap I'|=\ell-1$, let
$\epsilon(I,I')=(-1)^{p+q}$, where $K=I\cup I'$, $I=K_p$, and
$I'=K_q$.

\begin{thm}[\cite{AK}] \label{thm:gpGM}
For $\cG$--nonresonant weights $\la$, the Gauss-Manin connection on
the (trivial) vector bundle $\pi^\ell:\b{H}^\ell \to \sfB(\cG)$ has
connection $1$-form
\[
\sum_{J \subset [n+1]} d\log \sfD_J \otimes \Omega_\cG(\Gamma_{\!J}),
\]
where the sum is over all subsets $J$ of
$[n+1]$ of cardinality $\ell+1$ and $d \log \sfD_J$ denotes a $1$-form
on $(\CP^\ell)^n$ with a simple logarithmic pole along the divisor
$\sfD_J$.  The connection matrices $\Omega_\cG(\Gamma_{\!J})$, acting
on the $\beta$\textbf{nbc} basis $\{\eta_I \mid 2 \le i_1<\dots <
i_\ell \le n\}$ of $H^\ell(\sfE^{n,\ell}_{\b{y}};\LL_\b{y})$, are
given by the following formulas.

\smallskip

\noindent $\b{1.}$\quad If $J \cap \{1,n+1\}=\emptyset$, then
\begin{align*}
\Omega_\cG(\Gamma_{\!J})(\eta_I)&=
\begin{cases}
\sum_{p=1}^{\ell+1}\epsilon(I,J_p) \la_{J \setminus J_p} \eta_{J_p}
&\text{if $I\subset J$,}\\
0\phantom{12345678901234567890123456789012345678}
&\text{otherwise.}
\end{cases} \\
\intertext{$\b{2.}$\quad If $J =J'\cup\{n+1\}$, where
$J'=\{j_1,\dots,j_\ell\}$, $2\le j_1< \dots < j_\ell \le n$, then}
\Omega_\cG(\Gamma_{\!J})(\eta_I)&=
\begin{cases}
-\left(\sum_{j \in [n] \setminus I}\la_j\right) \eta_{I}
&\text{if $I= J'$,}\\
-\epsilon(I,J') \la_{I \setminus I\cap J'}\eta_{J'}
&\text{if $|I\cap J'|=\ell-1$,}\\
0\phantom{12345678901234567890123456789012345678}
&\text{otherwise.}
\end{cases}
\intertext{$\b{3.}$\quad If $J \cap\{1,n+1\}=\{1\}$, then}
\Omega_\cG(\Gamma_{\!J})(\eta_I)&=
\begin{cases}
\left(\sum_{j\in J}\la_j\right) \eta_{I} -
\sum_{|I\cap I'|=
\ell-1}\epsilon(I,I')\la_{I\setminus I \cap I'} \eta_{I'}
&\text{if $I=J_1$,}\\
0\phantom{12345678901234567890123456789012345678}
&\text{otherwise.}
\end{cases}
\intertext{$\b{4.}$\quad If $J =J'' \cup \{1,n+1\}$,
where $J''=\{j_2,\dots,j_\ell\}$,
$2 \le j_2 < \dots < j_\ell \le n$,
then}
\Omega_\cG(\Gamma_{\!J})(\eta_I)&=
\begin{cases}
-\la_{I\setminus I\cap J} \eta_{I} + \la_{I \setminus I\cap
J}\sum_{I\cap I'=I\cap J}\epsilon(I,I') \eta_{I'}
&\text{if $J'' \subset I$,}\\
0\phantom{12345678901234567890123456789012345678}
&\text{otherwise.}
\end{cases}
\end{align*}
\end{thm}

We combine these formulas with Theorem \ref{thm:lincomb}, Theorem
\ref{thm:nonresprojection}, and Theorem \ref{thm:GMCD} to obtain the
main result of this paper.  It is a combinatorial algorithm for
calculating Gauss-Manin connection matrices for an arbitrary
combinatorial type $\cT$, and $\cT$--nonresonant weights, in terms of
those for the type $\cG$ of general position arrangements.  Let $\cT$
be a combinatorial type which covers the type $\cT'$.  Let $\b{x}'$ be
a point in $\sfB(\cT')$, and $\c \in \pi_1(\sfB(\cT),\b{x})$ a simple
loop in $\sfB(\cT)$ about $\b{x}'$.  We agree to write
$\Omega_\cT(\cT')$ in place of $ \Omega_\cT(\c)$ in this situation.
Similarly, we write $ \Omega_\cG(J)$ in place of $
\Omega_\cG(\Gamma_{\!J})$.  Fix $\cT$--nonresonant (and hence
$\cG$--nonresonant) weights $\la$.  For $\b{x} \in \sfB(\cT)$ and
$\b{y} \in \sfB(\cG)$, by Theorems \ref{thm:vanish} and
\ref{thm:bnbc}, the corresponding local system cohomology groups of
the fibers, $H^q(\sfM_\b{x};\LL_\b{x})$ and
$H^q(\sfE^{n,\ell}_\b{y};\LL_\b{y})$, depend only on the combinatorial
types $\cT$ and $\cG$ and vanish for $q \neq \ell$.  So we write
$H^\ell(\sfM_\b{x};\LL_\b{x})=H^\ell(\cT)$ and
$H^\ell(\sfE^{n,\ell}_\b{y};\LL_\b{y})=H^\ell(\cG)$.  Denote the
$\beta$\textbf{nbc} bases of these local system cohomology groups by
$\beta$\textbf{nbc}$(\cT)$ and $\beta$\textbf{nbc}$(\cG)$
respectively.  This notation makes it clear that a Gauss-Manin
connection matrix depends only on the types involved and it is
obtained by a combinatorial formula.

\begin{thm} \label{thm:alg}
Let $\cT$ be a combinatorial type which covers the type $\cT'$.  Let
$\la$ be a collection of $\cT$--nonresonant weights, and let $P(\cT)$
be the matrix of the surjection $H^\ll(\cG) \twoheadrightarrow
H^\ll(\cT)$ in the respective $\beta$\textbf{nbc} bases.  Then a
Gauss-Manin connection matrix $\Omega_\cT(\cT')$ is determined by the
matrix equation
\[
P(\cT) \cdot \Omega_\cT(\cT') = \Bigl( \sum_{J \in \dep(\cT',\cT)}
m_J \cdot \Omega_\cG(J) \Bigr) \cdot P(\cT).
\]
\end{thm}
\noindent Illustrations of the algorithm provided by this result are
presented
in Section \ref{sec:examples}.

\begin{rem} \label{rem:hiro}
For an arbitrary combinatorial type $\cT$, and $\cT$--nonresonant
weights $\la$, Terao \cite{T1} shows that the Gauss-Manin connection
on the bundle $\pi^\ll:\b{H}^\ll\to\sfB(\cT)$ is determined by a
connection $1$-form $\sum d\log{\sf D}_j \otimes \Omega_j$, where
$\Omega_j \in \End_\C H^\ell(\sfM_\b{x};\LL_\b{x})$, each
$d\log\sfD_j$ denotes a $1$-form on $\overline\sfB(\cT)$ with a simple
logarithmic pole along the irreducible component ${\sf D}_j$ of the
codimension one divisor
$\sfD(\cT)=\overline\sfB(\cT)\setminus\sfB(\cT)$, and the sum is over
all such irreducible components.  Theorem \ref{thm:alg} computes the
Gauss-Manin connection matrices $\Omega_j$ here.  Since our work is
local, our results are independent of Terao's global theorem.
\end{rem}

\begin{rem}
Note that there are two floating inputs in Theorem \ref{thm:alg}: the
(STV)-condition, which ensures nonresonance, and the
$\beta$\textbf{nbc} basis, which makes the surjections explicit.
Should more relaxed conditions for nonresonance be discovered, and
suitable cohomology bases be constructed, where the surjection may be
made explicit, our methods will apply equally well in the new setting.

For an arbitrary local system $\LL$ of rank greater than one on the
complement of an arrangement of type $\cT$, the results presented here
need not apply, since such a local system is not in general induced by
a local system on the complement of a Boolean arrangement.  However,
if $\LL$ is an abelian local system of arbitrary rank, then $\LL$ is
induced by a local system on the complement of a Boolean arrangement.
The methods of this paper may be applied in this generality in those
instances where $\LL$ is nonresonant and a higher rank analogue of
Theorem \ref{thm:nonresprojection} holds.

In both of the above situations, analogues of the Aomoto-Kita formulas
for $\cG$--nonresonant weights, recorded in Theorem \ref{thm:gpGM}, would
be required to  express the Gauss-Manin connection matrices
for type $\cT$ in terms of those for type $\cG$ explicitly.
\end{rem}

\section{Examples} \label{sec:examples}
\subsection{Codimension one} \label{subsec:codim1}
The moduli space $\sfB(\cT)$ of a combinatorial type $\cT$ of
essential arrangements of $n$ hyperplanes in $\C^\ell$ is of
codimension one in $(\CP^\ell)^n$ if the cardinality of $\dep(\cT)$ is
one, $\dep(\cT)=\{K\}$, where $K=\{k_1,\dots,k_{\ell+1}\}$, $1\le k_1
< \dots < k_{\ell+1} \le n+1$.  For these types, the Gauss-Manin
connection was determined by Terao \cite{T1}.  We now sketch how these
results may be recovered using the algorithm of Theorem \ref{thm:alg}.

There are two cases to consider: $n+1 \notin K$ and $n+1 \in K$.  By
permuting hyperplanes, we may assume that
$K=[\ell+1]=\{1,\dots,\ell+1\}$ if $n+1 \notin K$, and that
$K=[n-\ell+1,n+1]=\{n-\ell+1,\dots,n+1\}$ if $n+1\in K$.  If $n+1
\notin K$, we have $\beta\text{\textbf{nbc}}(\cT)=
\beta\text{\textbf{nbc}}(\cG) \setminus \{\eta_F\}$, where
$F=[2,\ell+1]=\{2,\dots,\ell+1\}$, while if $n+1 \in K$,
$\beta\text{\textbf{nbc}}(\cT)= \beta\text{\textbf{nbc}}(\cG)
\setminus \{\eta_L\}$, where $L=[n-\ell+1,n]=\{n-\ell+1,\dots,n\}$.
Denote the projection $P(\cT):H^\ell(\cG) \twoheadrightarrow
H^\ell(\cT)$ of Theorem \ref{thm:alg} by $P_{F}$ if $n+1\notin K$ and
by $P_{L}$ if $n+1\in K$.  It is readily checked that these
projections are given~by
\begin{equation*} \label{eq:Aproj}
\begin{matrix}
P_F(\eta_I) = \eta_I \hfill & \text{if $I\neq F$,\quad} &
P_L(\eta_I)=\eta_I & \text{if $I \neq L$,} \\[2pt]
P_F(\eta_I)=\sum_{j=2}^{\ell+1} \sum_{q={\ell+2}}^n (-1)^{j+\ell}
\frac{\la_j}{\la_{[\ell+1]}} \eta_{F_{j-1}^q} & \text{if $I=F$,}\hfill
& P_L(\eta_I)=0 \hfill & \text{if $I=L$,}
\end{matrix}
\end{equation*}
where $F_{j-1}^q=F_{j-1}\cup\{q\}=\{2,\dots,\widehat j,
\dots,\ell+1,q\}$ for $\ell+2 \le q \le n$.

As shown by Terao \cite{T1}, the moduli space $\sfB(\cT)$ is smooth,
and the irreducible components of the divisor
$\sfD(\cT)=\overline{\sfB}(\cT) \setminus \sfB(\cT)$ correspond to the
combinatorial types $\cT'$ which are covered by $\cT$, and are
characterized as follows:
\begin{enumerate}
\item[(1)] $\dep(\cT',\cT) = \{J\}$, where $J\subset
[n+1]$, $|J|=\ell+1$, and $|K \cap J| \le \ell-1$
\item[(2)] $\dep(\cT',\cT) =\{K_p^q \mid q \in [n+1]
\setminus K\}$ for each $p$, $1\le p \le \ell+1$
\item[(3)] $\dep(\cT',\cT) = \{K_p^q \mid 1\le p \le
\ell+1\}$ for each $q \in [n+1] \setminus K$
\end{enumerate}
In types (2) and (3), $K_p^q=K_p \cup\{q\}=\{k_1,\dots,\widehat k_p,
\dots,k_{\ell+1},q\}$ for $1\le p \le \ell+1$ and $q \in [n+1]
\setminus K$.  If $\ell=1$, type (2) does not appear.

For each $J \in \dep(\cT',\cT)$, the order of vanishing of the
restriction of $\Delta_J$ to $\overline{\sfB}(\cT)$ along $\sfB(\cT')$
is equal to one, $m_J=1$, see \cite[\S10.4]{OT2}.  By Theorem
\ref{thm:alg}, a corresponding Gauss-Manin connection matrix
$\Omega_\cT(\cT') \in \End_\C(H^\ell(\cT))$ satisfies
\begin{equation} \label{eq:codim1GM}
P(\cT)\cdot \Omega_\cT(\cT') = \Bigl( \sum_{J \in \dep(\cT',\cT)}
\Omega_\cG(J)\Bigr) \cdot P(\cT),
\end{equation}
where $P(\cT)=P_F$ if $n+1\notin K$ and $P(\cT)=P_L$ if $n+1 \in K$.

It is an exercise in linear algebra to recover Terao's calculation of
the connection matrix $\Omega_\cT(\cT')$ from \cite{T1}.  We give
several small, explicit examples next, and leave this elementary,
albeit lengthy, exercise in the general case to the interested reader.

\begin{exm} \label{exm:123}
Let $\cT$ be the combinatorial type of the arrangement $\A$ of $4$
lines in $\C^2$ depicted in Figure \ref{fig:123}.
Here $\sfB(\cT)$ is codimension one in $(\CP^2)^4
=\overline{\sfB}(\cG)$.

\begin{figure}[h]
\setlength{\unitlength}{.45pt}
\begin{picture}(300,130)(-200,-110)
\put(-350,-100){\line(0,1){100}} \put(-420,-100){\line(1,1){100}}
\put(-280,-100){\line(-1,1){100}}\put(-420,-60){\line(1,0){140}}
\put(-320,5){3}\put(-280,-55){4}\put(-390,5){1}\put(-355,5){2}
\put(-357,-122){$\mathcal A$}

\put(-150,-100){\line(0,1){100}}\put(-220,-100){\line(1,1){100}}
\put(-80,-100){\line(-1,1){100}}\put(-185,-100){\line(1,1){100}}
\put(-120,5){3}\put(-80,5){4}\put(-190,5){1}\put(-155,5){2}
\put(-157,-122){${\mathcal A}_{1}$}

\put(250,-100){\line(0,1){100}}\put(180,-100){\line(1,1){100}}
\put(320,-100){\line(-1,1){100}}\put(180,-30){\line(1,0){140}}
\put(280,5){3}\put(320,-25){4}\put(210,5){1}\put(245,5){2}
\put(243,-122){${\mathcal A}_{3}$}

\put(50,-100){\line(0,1){100}}\put(-20,-100){\line(1,1){100}}
\put(-20,-60){\line(1,0){140}}
\put(80,5){3}\put(120,-55){4}\put(40,5){12}
\put(43,-122){${\mathcal A}_{2}$}
\end{picture}
\caption{A Codimension One Arrangement and Three Degenerations}
\label{fig:123}
\end{figure}
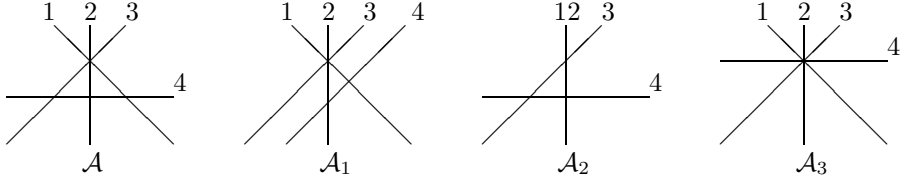

A collection of weights $\la$ is $\cT$--nonresonant if
$\la_j\ (1\le j \le 5),\ \la_1+\la_2+\la_3 \notin \Z_{\ge 0}$.
For such weights, $\beta$\textbf{nbc}$(\cT)=
\{\eta_{2,4},\eta_{3,4}\}$, and $\beta$\textbf{nbc}$(\cG)
=\{\eta_{2,3},\eta_{2,4},\eta_{3,4}\}$.  The projection
$P(\cT):H^2(\cG)\twoheadrightarrow H^2(\cT)$ has matrix
\[
P(\cT)=\begin{pmatrix}
\frac{-\la_3\phantom{-}}{\la_1+\la_2+\la_3} &
\frac{\la_2}{\la_1+\la_2+\la_3} \\
1 & 0 \\ 0 & 1 \end{pmatrix}.
\]

The arrangements $\A_{j}$ in Figure \ref{fig:123} represent types
$\cT'=\cT_{j}$ corresponding to irreducible components of the divisor
$\sfD(\cT)=\overline{\sfB}(\cT) \setminus \sfB(\cT)$.  The Gauss-Manin
connection matrices $\Omega_{\cT}(\cT_{j})$ are determined by Theorem
\ref{thm:alg} as follows.

\smallskip

\noindent{(1)}\quad
$\dep(\cT_{1},\cT)=\{345\}$, $m_{345}=1$,
$\Omega_\cG(345)= \smallmath{\begin{pmatrix} 0&0&-\la_2\\
0&0&\phantom{-}\la_2\\ 0&0&-\la_1-\la_2 \end{pmatrix}}$,\\[4pt]
\phantom{(1)}\quad $\implies \Omega_\cT(\cT_{1}) =
\begin{pmatrix} 0&\phantom{-}\la_2\\
0&-\la_1-\la_2 \end{pmatrix}$.\\[4pt]
(2)\quad $\dep(\cT_{2},\cT)=\{124,125\}$,
$m_{124}=m_{125}=1$,
$\Omega_\cG(125)=\smallmath{\begin{pmatrix}
-\la_3&-\la_3&0\\ -\la_4&-\la_4&0\\ \phantom{-}0&\phantom{-}0&0
\end{pmatrix}}$,\\[6pt]
\phantom{(2)}\quad
$\Omega_\cG(124)=
\smallmath{\begin{pmatrix} 0&0&0\\
\la_4&\la_1+\la_2+\la_4&\la_2\\ 0&0&0 \end{pmatrix}}
 \implies \Omega_\cT(\cT_{2}) = \begin{pmatrix}
\la_1+\la_2&\la_2\\ 0&0 \end{pmatrix}$.\\[6pt]
(3)\quad $\dep(\cT_{3},\cT)
=\{124,134,234\}$, $m_{124}=m_{134}=m_{234}=1$,\\[4pt]
\phantom{(3)}\quad
$\Omega_\cG(134)= \smallmath{\begin{pmatrix}
\phantom{-}0&0&0\\ \phantom{-}0&0&0\\
-\la_4&\la_3&\la_1+\la_3+\la_4
\end{pmatrix}}$,
$\Omega_\cG(234)=\smallmath{\begin{pmatrix}
\phantom{-}\la_4&-\la_3&\phantom{-}\la_2\\
-\la_4&\phantom{-}\la_3&-\la_2\\
\phantom{-}\la_4&-\la_3&\phantom{-}\la_2 \end{pmatrix}}$, \\[4pt]
\phantom{(3)}\quad
$\implies \Omega_\cT(\cT_{3}) = \begin{pmatrix}
\la_1+\la_2+\la_3+\la_4&0\\ 0&\la_1+\la_2+\la_3+\la_4
\end{pmatrix}$.
\end{exm}

\subsection{A Selberg arrangement} \label{subsec:selberg}
\begin{figure}[h]
\setlength{\unitlength}{.45pt}
\begin{picture}(300,130)(-200,-110)
\put(-200,-100){\line(0,1){110}}\put(-270,-100){\line(1,1){110}}
\put(-230,-100){\line(0,1){110}}\put(-270,-30){\line(1,0){110}}
\put(-270,-60){\line(1,0){110}}
\put(-155,15){5}\put(-155,-25){4}\put(-155,-55){3}\put(-235,15){1}
\put(-205,15){2}
\put(-224,-122){${\mathcal A}$}

\put(100,-100){\line(0,1){110}}
\put(70,-100){\line(0,1){110}}
\put(30,-60){\line(1,0){110}}
\put(140,-55){345}\put(65,15){1}\put(95,15){2}
\put(75,-122){${\mathcal A}'$}
\end{picture}
\caption{A Selberg Arrangement and One Degeneration}
\label{fig:selberg}
\end{figure}
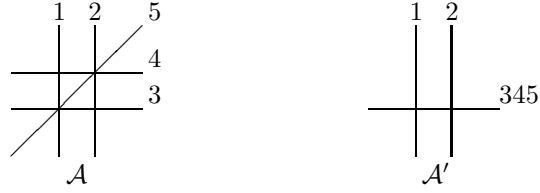

We conclude with an example of type $\cS$ for which $\sfB(\cS)$
 has higher codimension in  $\overline{\sfB}(\cG)$.
Let $\cS$ be the
combinatorial type of the Selberg arrangement $\A$ in $\C^2$ with
defining polynomial
$$
Q(\A)=u_1u_2(u_1-1)(u_2-1)(u_1-u_2)
$$
depicted in
Figure \ref{fig:selberg}.
See \cite{Ao,SV,JK} for detailed studies of the Gauss-Manin
connections arising in the context of Selberg arrangements.

A collection of weights $\la$ is $\cS$--nonresonant if
\[
\la_j\ (1\le j\le 6),\ \la_1+\la_2+\la_6,\ \la_1+\la_3+\la_5,
\ \la_2+\la_4+\la_5,\ \la_3+\la_4+\la_6\notin\Z_{\ge 0}.
\]
For such weights, $\beta$\textbf{nbc}$(\cS)
=\{\eta_{2,4},\eta_{2,5}\}$, $\beta$\textbf{nbc}$(\cG) = \{\eta_{i,j}
\mid 2 \le i<j\le 5\}$, and the projection $P(\cS):H^2(\cG)
\twoheadrightarrow H^2(\cS)$ is given by
\[
P(\cS)(\eta_{i,j})=
\begin{cases}
-\eta_{2,4}-\eta_{2,5} & \text{if $\{i,j\}=\{2,3\}$,}\\
\eta_{i,j} & \text{if $\{i,j\}=\{2,4\}$ or $\{2,5\}$,}\\
0 & \text{if $\{i,j\}=\{3,4\}$,}\\
\frac{(\la_3\la_5-\la_2\la_5) \eta_{2,4} -
(\la_2\la_3+\la_3\la_4+\la_2\la_5)\eta_{2,5}}{\la_2(\la_1+\la_3+\la_5)}
&\text{if $\{i,j\}=\{3,5\}$,}\\
\frac{-\la_5\eta_{2,4}+\la_4\eta_{2,5}}{\la_2}
&\text{if $\{i,j\}=\{4,5\}$.} \end{cases}
\]

The arrangement $\A'$ in Figure \ref{fig:selberg} represents one type
$\cS'$ covered by $\cS$.  Here
$\dep(\cS',\cS)=\{134,145,234,235,345,356,456\}$.  Since lines
$3$, $4$, and $5$ coincide in type $\cS'$, the order of vanishing of
the restriction of $\Delta_{345}$ to $\overline{\sfB}(\cS)$ along
$\sfB(\cS')$ is $m_{345}=2$.  For all other $J \in \dep(\cS',\cS)$, we
have $m_J=1$.  Theorem \ref{thm:alg} yields
\[
\Omega_\cS(\cS') = \begin{pmatrix} \la_3+\la_4+\la_5 & 0 \\
0 & \la_3+\la_4+\la_5 \end{pmatrix}.
\]

\bibliographystyle{amsalpha}

\end{document}